\documentclass[11pt]{article} % use larger type; default would be 10pt

\usepackage[utf8]{inputenc} % set input encoding (not needed with XeLaTeX)

\usepackage{mathtools}
\mathtoolsset{showonlyrefs}
\usepackage{geometry} % to change the page dimensions
\usepackage{amsthm, amsfonts, amssymb, amsmath,enumitem, natbib, bbm}
\newtheorem{thm}{Theorem}
\newtheorem{lem}{Lemma}

\newcounter{myalgctr}
\newenvironment{rem}{%      define a custom environment
   \vskip1mm\indent%         create a vertical offset to previous material
   \refstepcounter{myalgctr}% increment the environment's counter
   \textbf{Remark \themyalgctr}% or \textbf, \textit, ...
   }{\hfill$\diamond$\par}  %          create a vertical offset to following material
\numberwithin{myalgctr}{section}
\providecommand{\norm}[1]{\left\lVert#1\right\rVert}

\DeclareRobustCommand{\rchi}{{\mathpalette\irchi\relax}}
\newcommand{\irchi}[2]{\raisebox{\depth}{$#1\chi$}}
\renewcommand{\hat}{\widehat}
\usepackage{graphicx} % support the \includegraphics command and options
\usepackage{color, float}
\usepackage{booktabs} % for much better looking tables
\usepackage{array} % for better arrays (eg matrices) in maths
\usepackage{todonotes}
\usepackage{verbatim} % adds environment for commenting out blocks of text & for better verbatim
\usepackage{subfig} % make it possible to include more than one captioned figure/table in a single float
\usepackage[colorlinks=true, a4paper=true, pdfstartview=FitV,
linkcolor=blue, citecolor=blue, urlcolor=blue]{hyperref}
\makeatletter
\def\namedlabel#1#2{\begingroup
    #2%
    \def\@currentlabel{#2}%
    \phantomsection\label{#1}\endgroup
}
\makeatother
\newcommand{\vertiii}[1]{{\left\vert\kern-0.25ex\left\vert\kern-0.25ex\left\vert #1
    \right\vert\kern-0.25ex\right\vert\kern-0.25ex\right\vert}}
\usepackage[nottoc,notlof,notlot]{tocbibind} % Put the bibliography in the ToC
\usepackage[titles,subfigure]{tocloft} % Alter the style of the Table of Contents

\def\cred{\color{black}} %% NOTE (AC, 3/23/25): Introducing this color command for edits (typos, additions etc.). Feel free to use this or another color for your own edits.

\title{Tail Bounds for Canonical $U$-Statistics and $U$-Processes with Unbounded Kernels\thanks{An initial version of this work was available \href{https://faculty.wharton.upenn.edu/wp-content/uploads/2018/10/Chakrabortty-UStat-Draft.pdf}{here}.
This draft is a revised version with some edits.}%\emph{Note.} This is a working paper and the current draft is certainly not complete. The authors accept full responsibilities for any errors in this incomplete and unpublished manuscript\smallskip}
}%

\author{Abhishek Chakrabortty\thanks{Department of Statistics, Texas A\&M university. Email: \href{mailto:abhishek@stat.tamu.edu}{abhishek@stat.tamu.edu}}
 ~and~ Arun Kumar Kuchibhotla\thanks{Department of Statistics \& Data Science, Carnegie Mellon University. Email: \href{mailto:arunku@cmu.edu}{arunku@cmu.edu}}}

\date{\today} % Activate to display a given date or no date (if empty),
\begin{document}

\maketitle

\vspace{-0.1in}
\begin{abstract}
In this paper, we prove exponential tail bounds for canonical (or degenerate) $U$-statistics and $U$-processes under exponential-type tail assumptions on the kernels. Most of the existing results in the relevant literature often assume bounded kernels or obtain sub-optimal tail behavior under unbounded kernels. We obtain sharp rates and optimal tail behavior under sub-Weibull kernel functions. Some examples from nonparametric {\cred and semiparametric statistics} %nonparametric regression
literature are considered.
\end{abstract}

%%\par\smallskip

%\vspace{-0.075in}
\section{Introduction and Motivation}\label{sec:IntroMotiv}

%{\bf Convention:} {\cred edits in red}, {\cmag comments in magenta and with **} -- Abhishek (3/24/25)

In this paper, we study moment and tail bounds of second-order degenerate $U$-statistics and $U$-processes. Averages, the simplest function of a collection of random variables, are sums with each summand depending only on one element of the collection. On the other hand, $U$-statistics depend on tuples of elements in the collection. Formally, second-order $U$-statistics based on the collection of random variables $Z_1, \ldots, Z_n$ is of the form
\begin{equation}\label{eq:DefUStat}
{U}_n = \sum_{1\le i\neq j\le n} f_{i,j}(Z_i, Z_j),
\end{equation}
for some functions $\{f_{i,j}:\, 1\le i\neq j\le n\}$. In this paper, we consider $U$-statistics defined on independent but possible non-identically distributed random variables $Z_1, \ldots, Z_n$ defined on some measurable space.
{\cred $U$-statistics, in general, are ubiquitous in statistical applications including, e.g., goodness-of-fit tests, two-sample tests using kernel-based distances as well as independence testing via permutation tests; see \cite{kimstatistical2020} for an overview of this literature.}

{\cred We motivate our interest in $U$-statistics using a few prototypical examples.}
%% Changing the paragraph here not in the sentence before.
Suppose $X_1, \ldots, X_n$ are independent and identically distributed {\cred (i.i.d.)} realizations of a random vector $X\in\mathbb{R}^p$ with Lebesgue density $f$. Consider the problem of estimating the {\cred \it quadratic functional}
\begin{equation}\label{eq:IntegratedSquareFunctional}
\Gamma(f) := \int_{\mathbb{R}^p} f^2(x)dx = \mathbb{E}\left[f(X)\right].
\end{equation}
A natural estimator for this functional is given by
\begin{equation}\label{eq:EstimatorIntegratedSquare}
\widehat{\Gamma}(f) := \frac{1}{n(n-1)h_n^p}\sum_{1\le i\neq j\le n}K\left(\frac{X_i - X_j}{h_n}\right) = \frac{1}{n}\sum_{i=1}^n \widehat{f}^{(-i)}(X_i),
\end{equation}
where $h_n$ represents the bandwidth and $\hat{f}^{(-i)}(\cdot)$ represents the leave-one-out kernel density estimator:
\begin{equation}\label{eq:LeaveOneKernelEstimator}
\widehat{f}^{(-i)}(x) = \frac{1}{(n-1)h_n^p}\sum_{j = 1, j\neq i}^n K\left(\frac{X_j - x}{h_n}\right).
\end{equation}
Here the function $K(\cdot)$ is assumed to be symmetric and satisfies $\int_{\mathbb{R}^p} K(x)dx = 1$. This estimator was introduced by~\cite{HallMarron87} and was studied thoroughly (in terms of adaptivity) for $p = 1$ in~\cite{Gine08}.

Similarly, to estimate {\cred \it integrals} involving the conditional expectation function from {\cred i.i.d.} realizations $(X_1, Y_1), \ldots, (X_n, Y_n)$ of $(X, Y)$, the following $U$-statistics appears:
\begin{equation}\label{eq:EstimatorConditionalMean}
{U}_n^{\star} := \frac{1}{n(n-1)h_n^p}\sum_{1\le i\neq j\le n} Y_iK\left(\frac{X_i - X_j}{h_n}\right)Y_j.
\end{equation}
{\cred Aside from these prototypical examples, various other examples of such $U$-statistics are encountered in the literature on integral approximation involving kernel smoothing estimators \citep{NEW05, DEL16} and the semiparametric inference literature on quadratic and integral-type functionals \citep{robins2016asymptotic}. In the latter literature, $U$-statistics of this type -- especially in their degenerate form (see below for the definition) -- are fundamentally involved in the analysis of so-called {\it doubly robust estimators} of certain functionals encountered in missing data or causal inference problems \citep{robins1994estimation, bang2005doubly}, as well as in the literature on adaptive estimation of functionals based on so-called {\it higher order influence functions} \citep{robins2008higher,robins2017minimax,liu2021adaptive}.}

%{\cmag **Will need to make some edits/additions here -- between above and below part -- to add some stuff based on references mentioned later (at end of Section \ref{subsec:literature}). Also may need to create a subsection below when starting to talk about {\it degenerate} $U$-statistics in particular. (May also want to defer some of the literature/appln. discussions to that part, specific to degenerate $U$-statistics, especially the applications from semi-parametric inference (like the cases where it really shows up as degenerate $U$-statistics))**}

Apart from {\cred the} nonparametric {\cred and semiparametric} statistics literature, second order $U$-statistics also {\cred arise} %araise  %% TYPO HERE (Fixed, AC 3/23)
in relation to {\cred \it Hanson-Wright-type inequalities}. The classical Hanson-Wright inequality concerns tail bounds for the quadratic form $G^{\top}AG$ where $G$ is a standard multivariate normal random vector in $\mathbb{R}^n$ and $A\in\mathbb{R}^{n\times n}$ is a positive semi-definite matrix; see Theorem 3.1.9 of \cite{GINE16}. For further applications of Hanson-Wright inequalities, see~\cite{RudelsonVershynin13} and~\cite{Spokoiny13}, {\cred as well as the recent work of \cite{he2024sparse} on sparse random vectors.} Note that for any random vector $Y\in\mathbb{R}^n$ and matrix $A\in\mathbb{R}^{n\times n}$
\begin{equation}\label{eq:QuadraticFormToUStat}
Y^{\top}AY = \sum_{1\le i, j\leq n} Y_iA(i,j)Y_j,
\end{equation}
where $A(i,j)$ represents the $i$-th row, $j$-th column entry in the matrix $A.$

Motivated by the examples above, we study the properties of the $U$-statistic $U_n$. Before proceeding further, we briefly discuss degenerate and non-degenerate $U$-statistics. See \citet[Chapter 5]{SERF80} for more details. This discussion proves that for a precise understanding of the tail behavior of a $U$-statistics it suffices to consider degenerate $U$-statistics. In fact, most of the asymptotic normality results related to $U$-statistics are shown by proving asymptotic negligebility of the degenerate $U$-statistics compared to the linear statistic; see, for example,~\cite{CHEN17}. This paper is partly motivated by the cases where such asymptotic negligebility may not hold. For example, in the context of estimating $\mu^2$ based on IID observations $X_1, \ldots, X_n$ with mean $\mu$, the unbiased estimator $\binom{n}{2}^{-1}\sum_{i\neq j}X_iX_j$ exhibits a phase transition at $\mu = O(n^{-1/2})$ in terms of rate and also the limiting distribution.

\vspace{-0.05in}
\paragraph{Degenerate or Canonical $U$-statistics.} For any sequence of functions (called kernels) $f_{i,j}(\cdot, \cdot)$ and independent random variables $Z_1, \ldots, Z_n$, a $U$-statistic is given by
\[
T_n := \sum_{1\le i\neq j\le n}\,f_{i,j}(Z_i, Z_j).
\]
Note that the diagonal terms ($i = j$ cases) are ignored in the summation above. If these diagonal terms are included then the resulting statistic is called a $V$-statistic. The $U$-statistic $U_n$ is called {\cred \it degenerate} or {\cred \it canonical} %% Color only reflects italicizing here (no changes).
if the kernel functions satisfy
\begin{equation}\label{eq:DegenerateCondition}
\mathbb{E}\left[f_{i,j}(Z_i, Z_j)\big|Z_i\right] = \mathbb{E}\left[f_{i,j}(Z_i, Z_j)\big|Z_j\right] = 0,\quad\mbox{for all}\quad 1\le i\neq j\le n.
\end{equation}
If the kernel functions do not satisfy~\eqref{eq:DegenerateCondition}, then the corresponding $U$-statistic is called {\cred \it non-degenerate}.  %% Color only reflects italicizing here (no changes).
 It is not difficult to see that a non-degenerate $U$-statistic can be written as a sum of independent mean zero random variables and a degenerate $U$-statistic:
\begin{equation}\label{eq:NontoDegenerate}
T_n = \sum_{1\le i\neq j\le n} f_{i,j}^{D}(Z_i, Z_j) + \sum_{j = 1}^n g_j(Z_j) + \sum_{i=1}^n h_i(Z_i) =: \mathcal{U}_n(f) + T_n^{(1)} + T_n^{(2)},
\end{equation}
where
\begin{align}
% \begin{split}
f^D_{i,j}(Z_i, Z_j) &:= f_{i,j}(Z_i, Z_j) - \mathbb{E}\left[f_{i,j}(Z_i, Z_j)\big|Z_j\right] - \mathbb{E}\left[f_{i,j}(Z_i, Z_j)\big|Z_i\right] + \mathbb{E}\left[f_{i,j}(Z_i, Z_j)\right],\nonumber\\
g_j(Z_j) &:= \sum_{i = 1, i\neq j}^n \Big\{\mathbb{E}\left[f_{i,j}(Z_i, Z_j)\big|Z_j\right] - \mathbb{E}\left[f_{i,j}(Z_i, Z_j)\right]\Big\},\label{eq:SplitDegenerate}\\
h_i(Z_i) &:= \sum_{j = 1, j\neq i}^n \Big\{\mathbb{E}\left[f_{i,j}(Z_i, Z_j)\big|Z_i\right] - \mathbb{E}\left[f_{i,j}(Z_i, Z_j)\right]\Big\}.\nonumber
\end{align}
% \end{split}
% \end{equation}
It is clear from these expressions that the kernels $f_{i,j}^D(\cdot, \cdot)$ satisfy~\eqref{eq:DegenerateCondition} and so are degenerate kernels. Since $T_n^{(1)}$ and $T_n^{(2)}$ in~\eqref{eq:NontoDegenerate} are sums of independent random variables with mean zero, they can be understood easily from the classical results like the central limit theorem (asymptotically) and Bernstein/Hoeffding {\cred or more general} inequalities {\cred (non-asymptotically).} %(non-asumptotically). %% TYPO HERE (Fixed, AC 3/23)
For this reason, we focus mostly on the {\cred degenerate} %degerate %% TYPO HERE (Fixed, AC 3/23)
part of~\eqref{eq:NontoDegenerate} in the rest of the paper and derive non-asymptotic moment as well as tail bounds when the non-degenerate $U$-statistics is of the form~\eqref{eq:DefUStat}. Our main tool is the decoupling inequality proved in~\cite{DeLaPena92}. We refer to \citet[Chapter 3]{DeLaPena99} for more details regarding decoupling in $U$-statistics.

%% NOTE (AC, 3/28): Changing para here to highlight this is the next topic.
{\cred After} deriving non-asymptotic tail bounds for degenerate $U$-statistics, we provide the same for supremum of degenerate $U$-statistics over a function class. Suppose $\mathcal{F}_n$ is a class of sequence of functions (degenerate kernels) of type $f := \{f_{i,j}^{D}(\cdot, \cdot):\,1\le i\neq j\le n\}$ and define
\[
\mathcal{U}_n(f) := \sum_{1\le i\neq j\le n}\, f_{i,j}^D(Z_i, Z_j).
\]
Then $\{\mathcal{U}_n(f):\,f\in\mathcal{F}_n\}$ can be viewed as a process called the {\cred \it $U$-process} and we provide exponential tail bounds for the supremum:
\[
\mathcal{U}_n(\mathcal{F}) := \sup_{f\in\mathcal{F}_n}\,\left|\mathcal{U}_n(f)\right|.
\]
An important application would be the study of uniform-in-bandwidth properties of the estimator $\hat{\Gamma}(f)$ in~\eqref{eq:EstimatorIntegratedSquare}, that is,
\[
\sup_{h_n\in[a_n, b_n]}\left|\hat{\Gamma}(f; h_n) - \mathbb{E}\left[\hat{\Gamma}(f; h_n)\right]\right|,
\]
for some numbers $a_n, b_n\in(0, 1)$. Further applications can be found in~\citet[Section 5.5]{DeLaPena99} and~\cite{Major13}. As a final note, we mention that even though our techniques extend to $U$-statistics/processes of higher order, we restrict ourselves to second order $U$-statistics/processes {\cred for simplicity and ease of exposition}.

\subsection{\cred Related Literature}\label{subsec:literature} %{Literature Review} %%Changing subsection title (AC, 3/24)

%{\cmag **Changed subsection title from `Literature Review'**}

In this section, we review some of the by-now classical exponential tail bounds for degenerate $U$-statistics and supremum of $U$-processes. Proposition 2.3 of \cite{Arcones93} proved a Bernstein type inequality for degenerate $U$-statistics/processes. {\cred Specifically, for} the degenerate $U$-statistics
\[
U_n := n^{-1}\sum_{1\le i\neq j\le n}\, f(Z_i, Z_j),
\]
with {\cred i.i.d.} random variables $Z_1, \ldots, Z_n$, $\sigma^2 := \mathbb{E}f^2(Z_i, Z_j)$ and $\norm{f}_{\infty} \le C$, {\cred they show} there exists constants $c_1, c_2 > 0$ such that for any $t > 0$,
\[
\mathbb{P}\left(|U_n| \ge t\right) \le c_1\exp\left(-\frac{c_1t}{\sigma + (Ct^{1/2}n^{-1/2})^{2/3}}\right).
\]
This tail bound has two regimes: exponential and Weibull of order $2/3$. Because of the appearance of the variance, this tail bound provides the {\cred correct} %current
rate of convergence. Theorem 3.3 of \cite{Gine00} improved the tail bound by providing the optimal four regimes of the tail: Gaussian, exponential, Weibull of orders $2/3$ and $1/2$. \cite{HOUDRE03} gave an alternative proof to the result of \cite{Gine00} using martingale inequalities with explicit constants. In particular, Theorem 3.3 of \cite{Gine00} shows that for all $t\ge 0$,
\[
\mathbb{P}\left(|nU_n| \ge t\right) \le L\exp\left(-\frac{1}{L}\min\left\{\frac{t^2}{C^2}, \frac{t}{D}, \frac{t^{2/3}}{B^{2/3}}, \frac{t^{1/2}}{A^{1/2}}\right\}\right),
\]
for some constants $A, B, C, D$ and $L$. The main disadvantage of the results above is the restrictive boundedness assumption. Theorem 3.2 of \cite{Gine00} actually applies without the boundedness condition but the tail bound thus obtained is sub-optimal. For example, if $f(Z_i, Z_j) = Y_ig(X_i, X_j)Y_j$ where $Z_i = (X_i, Y_i)$, $\norm{g}_{\infty} \le C < \infty$ and $Y_i$'s are mean zero (conditionally) sub-Weibull variables of order $\alpha > 0$, that is, $\mathbb{P}(|Y_i| \ge t|X_i) \le 2\exp(-t^{\alpha})$. Then, Theorem 3.2 of \cite{Gine00} implies a tail bound of the form:
\[
\mathbb{P}\left(|nU_n| \ge t\right) \le L\exp\left(-\frac{1}{L}\min\left\{\frac{t^2}{C^2}, \frac{t}{D}, \frac{t^{\alpha_1}}{B^{\alpha_1}}, \frac{t^{\alpha_2}}{A^{\alpha_2}}\right\}\right),
\]
where $\alpha_1^{-1} = (3/2 + 1/\alpha)$ and $\alpha_2^{-1} = (2 + 2/\alpha)^{-1}$. This is sub-optimal in comparison with the results of \citet[Example 3]{Kole15}. On the other hand, the results of \cite{Kole15} do not get the correct rate of convergence as can be obtained from the results of \cite{Gine00}. This is because the bound of~\cite{Kole15} does not depend on the variance. We are not aware of any tail bounds in the literature that implies the correct rate of convergence as well as the optimal tail behavior.
{\cred 
We also note here the recent work of \cite{bakhshizadeh2023exponential} which appeared after the initial working version \citep{chakrabortty2018tail} of this preprint. While they do consider general unbounded kernels, their focus is primarily on exponential bounds and large deviation principles for non-degenerate $U$-statistics, different from ours.}

In regards to the tail bounds for {\cred degenerate} $U$-processes, some of the important works are \cite{Adamczak06}, \cite{CLEM08} and \cite{Major13}. The latter two papers only consider bounded kernels and the bounds of \cite{Adamczak06} are written in terms of functionals that are in general hard to control. The results of \cite{Major05} and \cite{Major13} apply only to bounded kernels and are written for VC {\cred classes} %class
 $\mathcal{F}_n$ but imply the correct rate of convergence. However, the results there do not show the optimal four regimes in the tail behavior. Theorem 11 of \cite{CLEM08} is written as a deviation inequality but does not imply the correct rate of convergence. For instance, if $f(X_i, X_j) = \varepsilon_i\varepsilon_jK((X_i - X_j)/h)$ with $\varepsilon_i$ being Rademacher random variables {\cred independent} %indepenedent %% TYPO HERE (Fixed, AC 3/23)
of $X_i\in\mathbb{R}^p$, then the rate of convergence of
\[
T_n := \sup_{h\in\{h_n\}}\,\left|\sum_{1\le i\neq j\le n}\varepsilon_i\varepsilon_jK\left(\frac{X_i - X_j}{h}\right)\right|,
\]
from Theorem 11 of \cite{CLEM08} is $n\norm{K}_{\infty} = O(n)$ (because of $Fn$ in the moment bound) but the correct rate of convergence is $nh_n^{p/2}$ (that can be obtained by calculating the variance). As in the case of $U$-statistics, we are not aware of any tail bound results that can obtain the correct rate of convergence and apply to unbounded kernels.
Using the techniques of truncation, decoupling technique and the entropy method of \cite{Bouch05}, we prove a deviation inequality for degenerate $U$-processes that implies the correct rate of convergence and the optimal tail behavior.

\vspace{-0.05in}
\paragraph*{\cred Organization.} The {\cred rest of the} %remaining
article is organized as follows. In Section~\ref{sec:UStat} we prove exponential tail bounds for second order degenerate $U$-statistics. In Section~\ref{sec:UProc} we prove a deviation bound for degenerate $U$-processes and also provide maximal inequalities to control the expectation of the maximum. %In Section~\ref{sec:Nonpar}, we apply our results to study {\color{red} an example} from non-parametric statistics. We conclude with some open questions and remarks in Section~\ref{sec:Conclusion}.
{\cred The proofs of all the results are distributed in Appendices \ref{sec:AppUStat}, \ref{appsec:UProc} and \ref{appsec:MaximalInequ}.}

\section[Tail Bounds for Degenerate U-Statistics]{Tail Bounds for Degenerate $U$--Statistics}\label{sec:UStat}
We prove two tail bounds for degenerate $U$-statistics. The first is a general result applicable to all kernels that are bounded above by a product kernel and the second result is for more structured kernels that are of importance in non- and semi-parametric estimation. Define a random variable $W$ to be {\cred \it sub-Weibull of order $\alpha > 0$} if $\norm{W}_{\psi_{\alpha}} < \infty,$ where $\psi_{\alpha}(x) = \exp(x^{\alpha}) - 1$ for $x \ge 0$ and
\[
\norm{W}_{\psi_{\alpha}} = \inf\left\{C \ge 0:\,\mathbb{E}\left[\psi_{\alpha}\left(|W|/C\right)\right] \le 1\right\}.
\]
Several properties of sum of independent sub-Weibull random variables are derived in~\cite{kuchibhotla2022moving}. %{\cred **Updated citation (final form)**}
The main focus of this section is to extend these results to degenerate $U$-statistics.

Consider a degenerate $U$-statistics
\[
U_n^D := \sum_{1\le i\neq j\le n} f_{i,j}(Z_i, Z_j),
\]
where $Z_1, \ldots, Z_n$ are independent random variables and $\{f_{i,j}(\cdot, \cdot):\, 1\le i \neq j\le n\}$ is a collection of degenerate (or canonical) kernels, i.e.,
\[
\mathbb{E}[f_{i,j}(Z_i, Z_j)|Z_i] = 0 = \mathbb{E}[f_{i,j}(Z_i, Z_j)|Z_j].
\]
We assume the following on the degenerate kernel $f_{i,j}$:
\begin{enumerate}[label=(\textbf{A\arabic*}),left=0cm]
\item For $1\le i\neq j\le n$, there exist non-negative functions $F_i(\cdot)$ and $G_j(\cdot)$ such that\label{assump:product-upper-bnd-kernel}
\[
|f_{i,j}(Z_i, Z_j)| \le F_i(Z_i)G_j(Z_j)\quad\mbox{and}\quad \|F_i(Z_i)\|_{\psi_{\alpha}} \le K_F,\; \|G_j(Z_j)\|_{\psi_{\beta}} \le K_G.
\]
\end{enumerate}
The first part of assumption~\ref{assump:product-upper-bnd-kernel} implies that the degenerate kernel $f_{i,j}$ can be expressed as $f_{i,j}(Z_i, Z_j) = F_i(Z_i)w_{i,j}(Z_i, Z_j)G_j(Z_j)$ for some collection of bounded kernels $\{w_{i,j}:1\le i\neq j\le n\}$. No further structure on $w_{i,j}$'s is required. In the second result we consider below, we place additional structure on $w_{i,j}$'s. Note that under assumption~\ref{assump:product-upper-bnd-kernel}, $F_i(\cdot)$ and $G_j(\cdot)$ are not uniquely defined because one can multiply $F_i(\cdot)$ by $2$ (say) and divide $G_j(\cdot)$ by $2$. In light of this, one can always assume $K_F = 1$; we do not assume this in the following. See~\ref{assump:weakend-product-kernel-assump} and the corresponding discussion.

The second part of assumption~\ref{assump:product-upper-bnd-kernel} means that
\[
\mathbb{E}\left[\exp\left(\frac{|F_i(Z_i)|^\alpha}{K_F^\alpha}\right)\right] \le 2\quad\mbox{and}\quad \mathbb{E}\left[\exp\left(\frac{|G_j(Z_j)|^\beta}{K_G^\beta}\right)\right] \le 2.
\]
Equivalently, $F_i(Z_i)$ is sub-Weibull$(\alpha)$ and $G_j(Z_j)$ is sub-Weibull$(\beta)$, in the terminology of~\cite{kuchibhotla2022moving}. To present the result, we define a few quantities. Let $(Z_1', Z_2', \ldots, Z_n')$ be an independent copy of $(Z_1, \ldots, Z_n)$.
\begin{equation}
\begin{split}
\Lambda_{1/2} &:= \left(\mathbb{E}\left[\sum_{1\le i\neq j\le n} f_{i,j}^2(Z_i, Z_j)\right]\right)^{1/2},\\
\Lambda_1 &= \|(f_{i,j})\|_{L^2\to L^2},\\
&:= \sup\left\{\mathbb{E}\sum_{1\le i\neq j\le n}f_{i,j}(Z_i, Z_j')\gamma_i(Z_i)\delta_j(Z_j'):\, \sum_{i=1}^n \mathbb{E}[\gamma_i^2(Z_i)] \le 1, \sum_{j=1}^n \mathbb{E}[\delta_j^2(Z_j')] \le 1\right\},\\
\Lambda_{\alpha} &:= \max_{1\le i\le n}\bigg\|\sum_{1\le j\le n, j\neq i} \mathbb{E}[f_{i,j}^2(Z_i, Z_j')|Z_i]\bigg\|_{\psi_{\alpha/2}}^{1/2},\\
\Lambda_{\beta} &:= \max_{1\le j\le n}\bigg\|\sum_{1\le i\le n, i\neq j}\mathbb{E}[f_{i,j}^2(Z_i, Z_j')|Z_j']\bigg\|_{\psi_{\beta/2}}^{1/2}.
\end{split}
\end{equation}
The quantities $\Lambda_{1/2}$ and $\Lambda_1$ also appear in the moment bound for degenerate $U$-statistics with bounded kernels; see Theorem 3.2 of~\cite{Gine00}.
Note that $\Lambda_{\alpha}$ can be trivially bounded as
\[
\Lambda_{\alpha} \le K_F\max_{1\le i\le n}\sup_{z}\left(\sum_{\substack{1\le j\le n,\\j\neq i}} \mathbb{E}\left[\frac{f_{i,j}^2(z, Z_j')}{F_i^2(z)}\right]\right)^{1/2}.
\]
Similar comment holds for $\Lambda_{\beta}$ as well. We use $\mathfrak{C}, \mathfrak{C}_{\alpha}, \mathfrak{C}_{\beta}$ and $\mathfrak{C}_{\alpha,\beta}$ to denote universal constants, constants depending on $\alpha, \beta, (\alpha,\beta)$, respectively. We now present the first main result.
\begin{thm}\label{thm:MainUStat-bounded-kernel}
For any collection of degenerate kernels $\{f_{i,j}:1\le i\neq j\le n\}$ satisfying assumption~\ref{assump:product-upper-bnd-kernel}, for every $p\ge1$,
\begin{equation}\label{eq:main-moment-inequality}
\begin{split}
&\left(\mathbb{E}\left[\left|\sum_{1\le i\neq j\le n} f_{i,j}(Z_i, Z_j)\right|^p\right]\right)^{1/p}\\ &\le \mathfrak{C}p^{1/2}\Lambda_{1/2} + \mathfrak{C}p\Lambda_1\\
&\quad+ \mathfrak{C}_{\beta}p^{1/2 + 1/\beta^*}(\log n)^{1/\beta}\Lambda_{\beta} + \mathfrak{C}_{\alpha}p^{1/2 + 1/\alpha^*}(\log n)^{1/2 + 1/\alpha}\Lambda_{\beta}\\
&\quad+ \mathfrak{C}_{\alpha,\beta}p^{1/\alpha^* + 1/\beta^*}K_FK_G(\log n)^{1/\alpha + 1/\beta + 1/\beta^*},
\end{split}
\end{equation}
where $\alpha^* = \min\{\alpha, 1\}$ and $\beta^* = \min\{\beta, 1\}$.
Consequently, for every $\delta\in[0,1]$, with probability at least $1-\delta$,
\begin{align*}
&\left|\sum_{1\le i\neq j\le n} f_{i,j}(Z_i, Z_j)\right|\\ &\le \mathfrak{C}(\log(1/\delta))^{1/2}\Lambda_{1/2} + \mathfrak{C}\log(1/\delta)\Lambda_1\\
&\quad+ \mathfrak{C}_{\beta}(\log(1/\delta))^{1/2 + 1/\beta^*}(\log n)^{1/\beta}\Lambda_{\beta} + \mathfrak{C}_{\alpha}(\log(1/\delta))^{1/2 + 1/\alpha^*}(\log n)^{1/2 + 1/\alpha}\Lambda_{\beta}\\
&\quad+ \mathfrak{C}_{\alpha,\beta}(\log(1/\delta))^{1/\alpha^* + 1/\beta^*}K_FK_G(\log n)^{1/\alpha + 1/\beta + 1/\beta^*}.
\end{align*}
\end{thm}
\begin{proof}
See %Section
Appendix~\ref{appsec:proof-thm:MainUStat-bounded-kernel} for a proof.
\end{proof}
Theorem~\ref{thm:MainUStat-bounded-kernel} reduces to Theorem 3.2 of~\cite{Gine00} by setting $\alpha = \beta = \infty$; note that if $\alpha = \beta = \infty$, then $\alpha^* = \beta^* = 1$ and the log factors in the result become 1. The proof of Theorem 3.2 of~\cite{Gine00} constitutes three main steps: (1) decoupling which implies that the study of $U_n^D$ follows from the study of the decoupled $U$-statistic
\[
\mathcal{U}_n := \sum_{1\le i\neq j\le n} f_{i,j}(Z_i, Z_j');
\]
(2) Rosenthal's or Pinelis' inequality conditional on $\{Z_1', \ldots, Z_n'\}$, and (3) an extension of Talagrand's inequality for supremum of averages (Proposition 3.1 of~\cite{Gine00}). The reason for sub-optimality of Theorem 3.2 of~\cite{Gine00} when applied to unbounded kernels stems from the sub-optimality of Rosenthal's inequality and Proposition 3.1 of~\cite{Gine00} for structured unbounded random variables. To elaborate Rosenthal's or Pinelis' inequality as stated in inequality (3.1) of~\cite{Gine00} is
\begin{equation}\label{eq:Rosenthal-Gine-3.1}
\mathbb{E}\left|\sum_{i=1}^n \xi_i\right|^p \le K^pp^{p/2}\left(\sum_{i=1}^n \mathbb{E}[\xi_i^2]\right)^{p/2} + K^pp^p\mathbb{E}\left[\max_{1\le i\le n}|\xi_i|^p\right],\quad\mbox{for}\quad p\ge2,
\end{equation}
for independent mean-zero random variables $\xi_1, \ldots, \xi_n$. If $|\xi_i| \le B$ almost surely, then this inequality is equivalent to Bernstein's inequality. However, for unbounded random variables $\xi_i$ (for examples, those satisfying only a sub-Gaussianity or sub-exponentiality), this inequality yields a sub-optimal moment bound. Our Theorem~\ref{thm:MainUStat-bounded-kernel} yields a sharper moment bound by replacing steps (2), (3) above with Theorem 3.4 and B.1 of~\cite{kuchibhotla2022moving}, respectively. Following the same logic, one can expand the scope of Theorem~\ref{thm:MainUStat-bounded-kernel} by replacing the sub-Weibull Orlicz norm with any other Orlicz norm for which Hoffmann-J {\o}rgensen inequality is satisfied. Recently,~\cite{adamczak2023orlicz} characterized the collection of Orlicz norms satisfying Hoffmann-J{\o}rgensen inequality which improves~\eqref{eq:Rosenthal-Gine-3.1} for structured unbounded random variables with the structure defined through an Orlicz norm. In particular, their result includes the classical sub-Weibull norm considered in Theorem~\ref{thm:MainUStat-bounded-kernel} and it also includes other Orlicz norms such as 
\[
\|W\|_{\Psi_{\alpha}} := \inf\{C > 0:\, \mathbb{E}[\Psi_{\alpha}(W/C)] \le 1\},
\]
where $\Psi_{\alpha}(x) = \exp(\log^{\alpha}(1 + x)) - 1$ for $x \ge 0$ and $\alpha \ge 1$. This Orlicz norm was introduced in~\cite{chamakh2021orlicz} and the random variables satisfying $\|W\|_{\Psi_{\alpha}} < \infty$ are called $\alpha$-heavy tailed. 

The result is {\em asymmetric} in $\alpha, \beta$ (even though assumption~\ref{assump:product-upper-bnd-kernel} is) because of the structure of the proof and in fact, the proof works with the following assumption (which is weaker than~\ref{assump:product-upper-bnd-kernel}):
\begin{enumerate}[label=(\textbf{A1$'$}),left=0cm]
\item There exist $\alpha, \beta > 0$ such that for $1\le i\neq j\le n$,\label{assump:weakend-product-kernel-assump}
\[
\widetilde{F}_i(z) = \max_{1\le j\le n}\|f_{i,j}(z, Z_j)\|_{\psi_{\beta}} < \infty,\quad\mbox{and}\quad \max_{1\le i\le n}\|\widetilde{F}_i(Z_i)\|_{\psi_{\alpha}} < \infty.
\]
\end{enumerate}
Under assumption~\ref{assump:weakend-product-kernel-assump}, the proof of Theorem~\ref{thm:MainUStat-bounded-kernel} goes through by taking $F_i(z) = \widetilde{F}_i(z),$ $K_F = \max_{1\le i\le n}\|\widetilde{F}_i(Z_i)\|_{\psi_{\alpha}}$, and $K_G = 1$; $G_j(\cdot)$ is never explicitly used in the proof. Unlike assumption~\ref{assump:product-upper-bnd-kernel}, assumption~\ref{assump:weakend-product-kernel-assump} is not symmetric in $\alpha$ and $\beta$ in the sense that assumption~\ref{assump:weakend-product-kernel-assump} does not imply the analogues assumption swapping the roles of $Z_i, Z_j$. More formally, assumption~\ref{assump:weakend-product-kernel-assump} is different from the following assumption:
\begin{enumerate}[label=(\textbf{A2$'$}),left=0cm]
\item There exist $\alpha, \beta > 0$ such that for $1\le i\neq j\le n$,\label{assump:weakend-product-kernel-assump-2}
\[
\widetilde{G}_j(z) = \max_{1\le i\le n}\|f_{i,j}(Z_i, z)\|_{\psi_{\alpha}} < \infty,\quad\mbox{and}\quad \max_{1\le i\le n}\|\widetilde{G}_j(Z_j)\|_{\psi_{\beta}} < \infty.
\]
\end{enumerate}
Note that assumptions~\ref{assump:weakend-product-kernel-assump} and~\ref{assump:weakend-product-kernel-assump-2} are significantly weaker than~\ref{assump:product-upper-bnd-kernel} as they do not require the product kernel upper bound on degenerate kernels. These modified assumptions allow application of results to non-decomposable kernels such as those in Lemma 2.19 of~\cite{bakhshizadeh2023exponential}. We mention in passing here that the product kernel assumption is not as significantly restrictive as it might seem because  any $U$-statistic can be decomposed into statistics satisfying~\ref{assump:product-upper-bnd-kernel} based on the Hilbert-Schmidt operator theory using the eigenfunctions of the kernel; see Section 5.5.2 of~\cite{SERF80} or Section 2 of~\cite{BentkusGotze1999} for details.

Under assumption~\ref{assump:product-upper-bnd-kernel}, one can obtain a symmetric result by applying Theorem~\ref{thm:MainUStat-bounded-kernel} by switching the roles of $\alpha, \beta$ and take the minimum of the two bounds. We do not present this for brevity. It is interesting to note that the tail exhibits five different behaviors including the commonly expected sub-Gaussian and sub-exponential tails. Because we did not make any assumption on the symmetry of the kernel, $\alpha$ and $\beta$ can be different. Under an assumption of symmetry, $\alpha = \beta$ and Theorem~\ref{thm:MainUStat-bounded-kernel} now yields a tail bound that only exhibits five regmies.

Assuming only~\ref{assump:product-upper-bnd-kernel}, Theorem~\ref{thm:MainUStat-bounded-kernel} provides a moment and tail bound for degenerate $U$-statistics. The appearance of the constants $\Lambda_{\alpha}$ and $\Lambda_{\beta}$ might make this result difficult to apply in some applications. For this reason, we provide our second result assuming a little more structure on the kernel.
Suppose we have $n$ independent random variables $Z_1 = (X_1, Y_1), Z_2 = (X_2, Y_2), \ldots, Z_n = (X_n, Y_n)$ on some measurable space and sequence of functions $\{w_{i,j}(\cdot, \cdot):\,1\le i \neq j\le n\}$. Consider, for functions $\phi(\cdot)$ and $\psi(\cdot)$, the $U$-statistic
\begin{equation}\label{eq:NonDegenerateUStat}
U_n := \sum_{1 \le i\neq j\le n} f_{i,j}(Z_i, Z_j),\quad\mbox{where}\quad f_{i,j}(Z_i, Z_j) := \phi(Z_i)w_{i,j}(X_i, X_j)\psi(Z_j).
\end{equation}
The kernels $f_{i,j}(\cdot, \cdot)$ are not required to be degenerate here. We will derive moment and tail bounds for the degenerate version of the $U$-statistics %$T_n^D$
{\cred $U_n^D$} given by
\[
U_n^D := \sum_{1\le i\neq j\le n} f_{i,j}^D(Z_i, Z_j),
\]
for the kernel $f_{i,j}^{D}(\cdot, \cdot)$ defined in~\eqref{eq:SplitDegenerate}. We first prove a basic lemma that reduces the problem of moment bounds on $U_n^{D}$ to a symmetrized version of $U_n$; see Theorem 3.5.3 of~\cite{DeLaPena99}. For any random variable $W$, set $\norm{W}_p = (\mathbb{E}[|W|^p])^{1/p}$ for $p\ge 1$.
\begin{lem}\label{lem:ReductionSymmetrization}
For any $p\ge 1$,
\[
\norm{U_n^D}_p \le C\norm{\sum_{1\le i\neq j\le n}\varepsilon_i\varepsilon_j'f_{i,j}(Z_i, Z_j')}_p,
\]
for Rademacher random variables $(\varepsilon_i, \varepsilon_i':\,1\le i\le n)$. Here $C$ can be taken to be $192$ and $Z_1' = (X_1', Y_1'), \ldots, Z_n' = (X_n', Y_n')$ represents an independent of $n$ independent random variables such that $Z_i$ is identically distributed as $Z_i$ for $1\le i\le n$.
\end{lem}
The proof of this lemma {\cred (given in Appendix \ref{appsec:proof-lem:ReductionSymmetrization})} %(in Section~\ref{sec:AppUStat}) 
is based on the by-now classical decoupling inequalities of \cite{DeLaPena92} and \citet[Chapter 3]{DeLaPena99}. The result also holds in case of degenerate $U$-processes and does not require the special structure of the kernels $f_{i,j}(\cdot, \cdot)$ in~\eqref{eq:NonDegenerateUStat}.

To prove moment and tail bounds for degenerate second order $U$-statistics with unbounded kernels, we use the following assumptions. Consider the following assumptions.
\begin{enumerate}[label = \bfseries(B\arabic*)]
\item There exists constants $0 < \alpha, \beta, C_{\phi}, C_{\psi} < \infty$ such that \label{assump:ConditionalTailAssumption}
\[
\max_{1\le i\le n}\,\mathbb{E}\left[\exp\left(\frac{|\phi(Z_i)|^{\alpha}}{C_{\phi}^{\alpha}}\right)\big|X_i\right] \le 2,\quad\mbox{and}\quad \max_{1\le i\le n}\,\mathbb{E}\left[\exp\left(\frac{|\psi(Z_i)|^{\beta}}{C_{\psi}^{\beta}}\right)\big|X_i\right] \le 2,
\]
hold almost surely.
\item The functions $\{w_{i,j}(\cdot, \cdot):\,1\le i\neq j\le n\}$ are all uniformly bounded, that is, \label{assump:BoundedKernels}
\[
\max_{1\le i\neq j\le n}\sup_{(x, x')\in\mathfrak{X}\times\mathfrak{X}}\,\left|w_{i,j}(x, x')\right| \le B_w.
\]
\end{enumerate}
The main technique in our proof is truncation and Hoffmann-J{\o}rgensen's inequality. Assumption~\ref{assump:ConditionalTailAssumption} implies that conditional on $X_i$'s the maximum of $\phi(Y_i)$ is at most a polynomial of $\log n$ (in rate). This along with Assumption~\ref{assump:BoundedKernels} allows us to apply truncation at this rate and study the truncated part using the sharp results of~\cite{Gine00}. The unbounded parts of smaller order are controlled using Hoffmann-J{\o}rgensen's inequality. The bound $B_w$ in Assumption~\ref{assump:BoundedKernels} is allowed to grown in $n$ and all the kernels are also allowed to be function of $n$. All the results to be presented here are non-asymptotic. 

Define
\begin{align*}
T_{\phi} := 8\mathbb{E}\left[\max_{1\le i\le n}\left|\phi(Z_i)\right|\big|X_1, \ldots, X_n\right],\quad T_{\psi} := 8\mathbb{E}\left[\max_{1\le i\le n}\left|\psi(Z_i)\right|\big|X_1, \ldots, X_n\right],
\end{align*}
and the truncated random variables
\begin{equation}\label{eq:TruncatedVariables}
\begin{split}
\Phi_{i,1} := \phi(Z_i)\mathbbm{1}\{|\phi(Z_i)| \le T_{\phi}\},\quad&\mbox{and}\quad \Phi_{i,2} := \phi(Z_i)\mathbbm{1}\{|\phi(Z_i)| > T_{\phi}\},\\
\Psi'_{j,1} := \psi(Z_j')\mathbbm{1}\{|\psi(Z_j')| \le T_{\psi}\},\quad&\mbox{and}\quad \Psi'_{j,2} := \psi(Z_j')\mathbbm{1}\{|\psi(Z_j')| > T_{\psi}\}.
\end{split}
\end{equation}
It is clear that $\phi(Z_i) = \Phi_{i,1} + \Phi_{i,2}$ and $\psi(Z_j') = \Psi'_{j,1} + \Psi'_{j,2}.$ Based on these, note that
\begin{equation}\label{eq:BasicDecomposition}
\begin{split}
\phi(Z_i)w_{i,j}(X_i, X_j')\psi(Z_j') &= \Phi_{i,1}w_{i,j}(X_i, X_j)\Psi'_{j,1} + \Phi_{i,2}w_{i,j}(X_i, X_j)\Psi'_{j,1}\\ &\qquad+
\Phi_{i,1}w_{i,j}(X_i, X_j)\Psi'_{j,2} +
\Phi_{i,2}w_{i,j}(X_i, X_j)\Psi'_{j,2}.
\end{split}
\end{equation}
The first term on the right hand side is bounded by $T_{\phi}B_wT_{\psi}$. The second and third terms are non-zero only when $\Phi_{i,2}$ and $\Psi'_{j,2}$, are respectively non-zero, which can only happen with only a small probability under Assumption~\ref{assump:ConditionalTailAssumption}. Finally, the fourth term can be non-zero only if both $\Phi_{i,2}$ and $\Psi'_{j,2}$ are non-zero which can {\cred happen} %happend %% TYPO HERE (Fixed, AC 3/23)
with even smaller probability. These four terms leads to four different degenerate $U$-statistics that will be controlled separately in Section~\ref{subsec:AppSecMainUStat} to prove the following result. We need the following notation: for $1\le i, j\le n$,
\[
\sigma_{i,\phi}^2(x) = \mathbb{E}[\phi^2(Z_i)\big|X_i = x]\quad\mbox{and}\quad\sigma_{j,\psi}^2(x) = \mathbb{E}[\psi^2(Z_j)\big|X_j = x].
\]
Define $\Lambda_2 := C_{\phi}C_{\psi}B_w(\log n)^{\alpha^{-1} + \beta^{-1}}$ and
\begin{align*}
\Lambda_{{1}/{2}} &:= \left(\sum_{1\le i\neq j\le n}\mathbb{E}\left[\sigma_{i,\phi}^2(X_i)w_{i,j}^2(X_i, X_j)\sigma_{j,\psi}^2(X_j)\right]\right)^{1/2},\\
\Lambda_1 &:= \sup\left\{\sum_{1\le i\neq j\le n}\mathbb{E}\left[q_i(X_i)\sigma_{i,\phi}(X_i)w_{i,j}(X_i, X_j)\sigma_{j,\psi}(X_j)p_j(X_j)\right]:\right.\\
&\qquad\qquad\quad\left.\sum_{j = 1}^n \mathbb{E}\left[q_i^2(X_i)\right] \le 1,\,\sum_{i=1}^n \mathbb{E}\left[p_j^2(X_j)\right] \le 1\right\},\\
\Lambda_{3/2}^{(\alpha)} &:= C_{\phi}(\log n)^{1/\alpha}\sup_{x} \max_{1\le i\le n} \left(\sum_{j = 1}^n \mathbb{E}\left[w_{i,j}^2(x, X_j)\sigma^2_{j,\psi}(X_j)\right]\right)^{1/2},\\
\Lambda_{3/2}^{(\beta)} &:= C_{\psi}(\log n)^{1/\beta}\sup_{x} \max_{1\le j \le n} \left(\sum_{i = 1}^n \mathbb{E}\left[w_{i,j}^2(X_i, x)\sigma^2_{i, \phi}(X_i)\right]\right)^{1/2},\\
\Lambda_{\alpha^*} &:= (\log n)^{1/2}\Lambda_{3/2}^{(\alpha)} + (\log n)\Lambda_2, \quad \mbox{\cred and}\\
\Lambda_{\beta^*} &:= (\log n)^{1/2}\Lambda_{3/2}^{(\beta)} + (\log n)\Lambda_2.
\end{align*}
The quantities $\Lambda_{1/2}, \Lambda_{1}, \Lambda_{3/2}^{(\alpha)}, \Lambda_{3/2}^{(\beta)}, \Lambda_2$ also appear in the case of bounded kernels as shown in Theorem 3.2 of~\cite{Gine00}.
\begin{thm}\label{thm:MainUStat}
Under Assumptions~\ref{assump:ConditionalTailAssumption} and~\ref{assump:BoundedKernels}, there exists constant $K > 0$ (depending only on $\alpha, \beta$) such that for all $p\ge1$
\begin{align*}
\norm{U_n^D}_p &\le K p^{1/2}\Lambda_{1/2} + K p\Lambda_1 + K p^{1/\alpha^*}\Lambda_{\alpha^*} + Kp^{1/\beta^*}\Lambda_{\beta^*}\\
&\qquad + K p^{1/2 + 1/\alpha^*}\Lambda_{3/2}^{(\alpha)} + K p^{1/2 + 1/\beta^*}\Lambda_{3/2}^{(\beta)} + Kp^{1/\alpha^* + 1/\beta^*}\Lambda_2.
\end{align*}
Here $\alpha^* := \min\{\alpha, 1\}$ and $\beta^* := \min\{\beta, 1\}$. By Markov's inequality, there exists a constant $K' > 0$ such that for any $t\ge 0$,
\begin{equation}\label{eq:TailBoundDegenerateStat}
\mathbb{P}\left(|U_n^D| \ge K'\mathcal{T}_{\alpha, \beta}(t)\right) \le 2\exp(-t),
\end{equation}
where
\begin{align*}
\mathcal{T}_{\alpha, \beta}(t) &:= \sqrt{t}\Lambda_{1/2} + t\Lambda_1 + t^{1/\alpha^*}\Lambda_{\alpha^*} + t^{1/\beta^*}\Lambda_{\beta^*}\\
&\qquad+ t^{1/2 + 1/\alpha^*}\Lambda_{3/2}^{(\alpha)} + t^{1/2 + \beta^*}\Lambda_{3/2}^{(\beta)} + t^{1/\alpha^* + 1/\beta^*}\Lambda_2.
\end{align*}
\end{thm}

\begin{proof}
{\cred See Appendix~\ref{subsec:AppSecMainUStat} for a proof.}
\end{proof}

\begin{rem}\,(Comparison with previous results)
As noted in the introduction, an important feature of our result is that the kernel is allowed to be unbounded with proper tail behavior. The tail of the degenerate $U$-statistics as shown in~\eqref{eq:TailBoundDegenerateStat} has \emph{seven} different regimes, the prominent ones being the Gaussian and exponential parts. These seven regimes collapse to \emph{five} if $\alpha = \beta$. In particular, if $\alpha = \beta \le 1$, then for $p\ge 1$,
\begin{align*}
\norm{U_n^D}_p &\le K \mathbf{p^{1/2}}\Lambda_{1/2} + K \mathbf{p}\Lambda_1 + K \mathbf{p^{1/\alpha}}\left[(\log n)^{1/2}\left\{\Lambda_{3/2}^{(\alpha)}  + \Lambda_{3/2}^{(\beta)}\right\} + (\log n)\Lambda_2\right]\\
&\quad + K \mathbf{p^{1/2 + 1/\alpha}}\left[\Lambda_{3/2}^{(\alpha)} + \Lambda_{3/2}^{(\beta)}\right] + K\mathbf{p^{1/\alpha + 1/\beta}}\Lambda_2.
\end{align*}
If $\alpha = \beta = \infty$, then our assumption~\ref{assump:ConditionalTailAssumption} implies boundedness of the kernels. In this case, only four regimes remain and these four regimes coincide with those shown in Theorem 3.2 of~\cite{Gine00}. Additionally in the case of bounded kernels ($\alpha = \beta = \infty$), Theorem~\ref{thm:MainUStat} essentially coincides with Theorem 3.2 of~\cite{Gine00} except for the additional $\sqrt{\log n}$ and $\log n$ factors. We believe these to be artifacts of our proof and closely following the proof of Theorem~\ref{thm:MainUStat-bounded-kernel}, they could be avoided.
\end{rem}

\section[Tail Bounds for Degenerate U-Processes]{Tail Bounds for Degenerate $U$--Processes}\label{sec:UProc}

In this section, we generalize Theorem~\ref{thm:MainUStat} to degenerate $U$-processes. Consider
\[
\mathcal{U}_{n}(\mathcal{W}) := \sup_{w\in\mathcal{W}}\left|\mathcal{U}_n(w)\right|,\quad\mbox{where}\quad \mathcal{U}_n(w) := \sum_{1\le i\neq j\le n} \varepsilon_i\phi(Z_i)w_{i,j}(X_i, X_j')\psi(Z_j')\varepsilon_j',
\]
for some function class $\mathcal{W}$ with elements of the type $w = (w_{i,j})_{1\le i\neq j\le n}$. If $\mathcal{W}$ is a singleton, then this reduces to the $U$-statistic studied in Section~\ref{sec:UStat}. Here $\varepsilon_1, \ldots, \varepsilon_n$ denote an independent sequence of Rademacher random variables as before. For simplicity, we consider the symmetrized version and by Lemma~\ref{lem:ReductionSymmetrization} the results also hold for the original degenerate $U$-process; see Theorem 3.5.3 of~\cite{DeLaPena99} for details.

$U$-processes were introduced in~\cite{Nolan87} to study cross-validation in the context of kernel density estimation. They studied uniform almost sure limit theorems and established the rate of convergence. These results parallel the Glivenko-Cantelli theorems well-known for empirical processes. Functional limit theorems were established in~\cite{nolan1988functional}. Exponential tail bounds that parallel the classical Bernstein's inequality for non-degenerate and degenerate $U$-statistics were given in~\cite{Arcones93}. They also established LLN and CLT type results under various metric entropy conditions. Most of these results require boundedness of the kernel functions. Being asymptotic in nature, some of these results can be extended to the case of unbounded kernels using a truncation argument. Finite sample concentration inequalities for degenerate unbounded $U$-processes are not readily available.

The only work (we are aware of) that provides general results for $U$-processes applicable to $\mathcal{U}_n(\mathcal{W})$ is~\cite{Adamczak06}. In this work, degenerate $U$-processes of arbitrary order were considered. However, the moment bounds for $U$-processes in this work depend further on the moment bounds of some complicated degenerate $U$-processes of lower order. Furthermore, the tail behavior thus obtained is not sharp for unbounded $U$-processes.

To avoid measurability issues for $\mathcal{U}_n(\mathcal{W})$, we use either of the following conventions. One simple assumption on $\mathcal{W}$ used in~\cite{VdvW96} that implies measurability is separability and in this case we can take $\mathcal{W}$ to be a dense countable subset of $\mathcal{W}$. Another convention used in~\cite{Tala14} is to define for any $\mathcal{W}$ and increasing function $f(\cdot)$,
\[
\mathbb{E}\left[f\left(\mathcal{U}_n(\mathcal{W})\right)\right] := \sup\left\{\mathbb{E}\left[f(\mathcal{U}_n(\mathcal{F}))\right]:\,\mathcal{F}\subseteq \mathcal{W}\mbox{ a finite subset}\right\}.
\]
Based on either convention, we treat $\mathcal{W}$ as a countable set for the remaining part of this section.

One ``simple'' way to obtain tail bounds for $\mathcal{U}_n(\mathcal{W})$ is via generic chaining as follows: First apply Theorem~\ref{thm:MainUStat} for $\mathcal{U}_n(w) - \mathcal{U}_n(w')$ for functions $w, w'\in\mathcal{W}$. The tail bound~\eqref{eq:TailBoundDegenerateStat} provides a mixed tail in terms of various semi-metrics on $\mathcal{W}$. Using these and following the proof of classical generic chaining bound (e.g., Theorem 3.5 of~\cite{Dirk15}), one can obtain tail bounds for $U$-processes in terms of $\gamma$-functionals; see~\cite{Tala14} and~\cite{Dirk15} for details. A problem with this approach is the complication in controlling the $\gamma$-functionals. This approach with Dudley's chaining (instead of generic chaining) was used for bounded kernel $U$-processes in~\cite{Nolan87} and~\cite{nolan1988functional}.

In the following, we first provide a deviation inequality for $\mathcal{U}_n(\mathcal{W})$ and then prove a maximal inequality to control the expectations appearing in the deviation inequality. For these results, we consider the following generalization of assumption~\ref{assump:BoundedKernels}.
\begin{enumerate}[label = \bfseries(B2$'$)]
\item The functions $\{w:\,w\in\mathcal{W}\}$ are all uniformly bounded, that is, \label{assump:BoundedKernelsUniform}
\[
\sup_{w\in\mathcal{W}}\sup_{(x, x')\in\mathfrak{X}\times\mathfrak{X}}\,\max_{1\le i\neq j\le n}\left|w_{i,j}(x, x')\right| \le B_{\mathcal{W}}.
\]
\end{enumerate}
We will use the notation of $\Phi_{i,1}, \Phi_{i,2}, \Psi_{i,1}', \Psi_{i,2}'$ given in~\eqref{eq:TruncatedVariables}. For the main result of this section, define
\begin{align*}
% E_{n,0}(\mathcal{W}) &:= \mathbb{E}\left[\sup_{w\in\mathcal{W}}\left|\sum_{1\le i\neq j\le n}\varepsilon_i\Phi_{i,1}w(X_i, X_j')\Psi_{j,1}'\varepsilon_j'\right|\right],\\
\Lambda_2(\mathcal{W}) &:= (\log n)^{\alpha^{-1} + \beta^{-1}}C_{\phi}C_{\psi}B_{\mathcal{W}},\\
E_{n,1}(\mathcal{W}) &:= C_{\psi}(\log n)^{1/\beta}\sup_{x\in\mathfrak{X}}\max_{1\le j\le n}\mathbb{E}\left[\sup_{w\in\mathcal{W}}\left|\sum_{i=1,i\neq j}^n \varepsilon_i\Phi_{i,1}w_{i,j}(X_i, x)\right|\right],\\
{E}_{n,2}(\mathcal{W}) &:= C_{\phi}(\log n)^{1/\alpha}\sup_{x\in\mathfrak{X}}\max_{1\le i\le n}\mathbb{E}\left[\sup_{w\in\mathcal{W}}\left|\sum_{j = 1, j\neq i}^n \varepsilon_j'\Psi_{j,1}'w_{i,j}(x, X_j')\right|\right],\\
\mathfrak{W}_{n,1}(\mathcal{W}) &:= \mathbb{E}\left[\sup_{w\in\mathcal{W}}\sup_{\{p_j\}}\sum_{1\le i\neq j\le n} \varepsilon_i\Phi_{i,1}\int p_j(x)\sigma_{j,\psi}(x)w_{i,j}(X_i, x)P_{X_j}(dx)\right],\\
\mathfrak{W}_{n,2}(\mathcal{W}) &:= \mathbb{E}\left[\sup_{w\in\mathcal{W}}\sup_{\{q_i\}}\sum_{1\le i\neq j\le n} \varepsilon_j'\Psi_{j,1}'\int q_i(x)\sigma_{i,\phi}(x)w_{i,j}(x, X_j')P_{X_i}(dx)\right],\\
%\end{align*}  %% NOTE (AC, 3/29): Breaking up the align environment here or else formatting looks bad.
%\begin{align*}
{\cred \Sigma_{n,1}^{1/2}(\mathcal{W})} &:= C_{\psi}(\log n)^{1/\beta}\sup_{x\in\mathfrak{X}}\sup_{w\in\mathcal{W}}\,\max_{1\le j\le n}\left(\sum_{i=1,i\neq j}^n \mathbb{E}[\sigma_{i,\phi}^2(X_i)w^2_{i,j}(X_i, x)]\right)^{1/2},\\
{\Sigma}_{n,2}^{1/2}(\mathcal{W}) &:= C_{\phi}(\log n)^{1/\alpha}\sup_{x\in\mathfrak{X}}\sup_{w\in\mathcal{W}}\max_{1\le i\le n}\left(\sum_{j = 1, j\neq i}^n \mathbb{E}\left[\sigma_{j,\psi}^2(X_j)w^2_{i,j}(x, X_j)\right]\right)^{1/2},\\
\norm{(\phi w\psi)_{\mathcal{W}}}_{2\to2} &:= \sup_{w\in\mathcal{W}}\sup_{\{q_i\}}\sup_{\{p_j\}}\sum_{1\le i\neq j\le n}\mathbb{E}\left[q_i(X_i)\sigma_{i,\phi}(X_i)w_{i,j}(X_i, X_j')\sigma_{j,\psi}(X_j')p_j(X_j')\right].
\end{align*}
Here in the definitions, the supremum over $\{q_i\}$ (or $\{p_j\}$) represents supremum over all function $(q_1, \ldots, q_n)$ (or $(p_1, \ldots, p_n)$)satisfying
\[
\sum_{i=1}^n \int q_i^2(x)P_{X_i}(dx) \le 1,\quad\mbox{and}\quad \sum_{j = 1}^n \int p_j^2(x)P_{X_i}(dx) \le 1,
\]
where $P_{X_i}(\cdot)$ denotes the probability measure of $X_i$. Note that $\norm{(\phi w\psi)_{\mathcal{W}}}_{2\to 2}$ is similar to $\Lambda_{1}$ defined for Theorem~\ref{thm:MainUStat}.
\begin{thm}\label{thm:MainUProc}
Under assumptions~\ref{assump:ConditionalTailAssumption} and~\ref{assump:BoundedKernelsUniform}, there exists a constant $K > 0$ (depending only on $\alpha, \beta$) such that for all $p\ge 1$
\begin{align*}
\norm{\mathcal{U}_n(\mathcal{W})}_p &\le K\mathbb{E}\left[\mathcal{U}_n^{(1)}(\mathcal{W})\right] + Kp^{1/2}(\mathfrak{W}_{n,1}(\mathcal{W}) + \mathfrak{W}_{n,2}(\mathcal{W})) + Kp\norm{(\phi w\psi)_{\mathcal{W}}}_{2\to2}\\
 % + Kp^{3/2}\left(\Sigma_{n,1}^{1/2}(\mathcal{W}) + \Sigma_{n,2}^{1/2}(\mathcal{W})\right) + Kp^2\Lambda_2(\mathcal{W})\\
&\quad+ Kp^{1/\alpha^*}\left[E_{n,2}(\mathcal{W}) + \Sigma_{n,2}^{1/2}(\mathcal{W})\sqrt{\log n} + \Lambda_2(\mathcal{W})\log n\right]\\
&\quad+ Kp^{1/\beta^*}\left[E_{n,1}(\mathcal{W}) + \Sigma_{n,1}^{1/2}(\mathcal{W})\sqrt{\log n} + \Lambda_2(\mathcal{W})\log n\right]\\
&\quad+ Kp^{1/2 + 1/\alpha^*}\Sigma_{n,2}^{1/2}(\mathcal{W}) + Kp^{1/2 + 1/\beta^*}\Sigma_{n,1}^{1/2}(\mathcal{W}) + Kp^{1/\alpha^* + 1/\beta^*}\Lambda_2(\mathcal{W}).
\end{align*}
\end{thm}
\begin{proof}
See Appendix~\ref{appsec:MainUProc} for a proof.
\end{proof}
If $\mathcal{W}$ is a singleton set, then the above result reduces to Theorem~\ref{thm:MainUStat}. From the moment bound above, it is easy to derive a tail bound using Markov's inequality. In comparison, we again get \emph{seven} different tail regimes that again reduce to \emph{five} if $\alpha = \beta$. Unlike the result of~\cite{Adamczak06}, the moment bound in Theorem~\ref{thm:MainUProc} only depends on some expectations. An additional advantage of Theorem~\ref{thm:MainUProc} is that all the expectations only involve bounded random variables.

\subsection{Maximal Inequality for Bounded Degenerate \textit{U}-Processes}\label{sec:MaxIneq}

To apply Theorem~\ref{thm:MainUProc}, we need to control various expectations appearing on the right hand side of the moment bound there. Expect for $\mathbb{E}[\mathcal{U}_n^{(1)}(\mathcal{W})]$, all the other quantities are maximal inequalities related to empirical processes. See~\cite{vdV11} and Lemmas 3.4.2-3.4.3 of~\cite{VdvW96} for maximal inequalities of empirical processes. In this section, we provide a maximal inequality for $\mathcal{U}_n^{(1)}(\mathcal{W})$. For independent and identically distributed random variables, \citet[Theorem 5.1]{CHEN17} %{\cred **Updated citation (journal form)}
provide a maximal inequality for degenerate $U$-processes of arbitrary order. This result is similar to Theorem 2.1 of~\cite{vdV11} for empirical processes. The same proof as in~\cite{CHEN17} does not provide the ``correct'' bound in the case of possibly non-identically distributed observations since they use Hoeffding {\cred averaging} %averging %% TYPO HERE (Fixed, AC 3/23)
which can lead to sub-optimal rate if the observations are not identically distributed. A modification of the proof leads to the maximal inequality below.

For any $\eta > 0$, function class $\mathcal{F}$ containing functions $f = (f_{i,j})_{1\le i\neq j\le n}:\rchi\times\rchi\to\mathbb{R}$ and a discrete probability measure $Q$ with support $\{z_1, \ldots, z_t\}$, let $N(\eta, \mathcal{F}, \norm{\cdot}_{2,Q})$ denotes the minimum $m$ such that there exists $f^{(1)}, f^{(2)}, \ldots, f^{(m)}\in\mathcal{F}$ satisfying
\[
\sup_{f\in\mathcal{F}}\inf_{1\le j\le m}\norm{f - f^{(j)}}_{2,Q} \le \eta,
\]
where for $f\in\mathcal{F}$,
\[
\norm{f}_{2,Q}^2 := \frac{\sum_{1\le i\neq j\le t} f_{i,j}^2(z_i, z_j)Q(\{z_i\})Q(\{z_j\})}{\sum_{1\le i\neq j\le t}Q(\{z_i\})Q(\{z_j\})}.
\]
Note that the right hand side is expectation with respect to the measure induced on $\{(z_i, z_j):1\le i\neq j\le t\}$. Define the uniform entropy integral needed for $U$-processes is given by
\[
J_2(\delta, \mathcal{F}, \norm{\cdot}_2) := \sup_{Q}\int_0^{\delta} \log N(\eta\norm{F}_{2,Q}, \mathcal{F}, \norm{\cdot}_{2,Q})d\eta.
\]
Here $F = (F_{i,j})_{1\le i\neq j\le n}$ represents the envelope function for $\mathcal{F}$ satisfying $|f_{i,j}(x, x')| \le F_{i,j}(x,x')$ for all $f\in\mathcal{F}, x,x'\in\rchi$ and the supremum is taken over all discrete probability measures $Q$ supported on $\rchi\times\rchi$.

The following Lemma proves a maximal inequality using Theorem 5.1.4 of \cite{DeLaPena99}. The proof is very similar to that of Theorem 5.1 of \cite{CHEN17} which itself was based on the proof of Theorem 2.1 of \cite{vdV11}.

\begin{thm}\label{thm:UniformEntropyUProcesses}
Suppose $\mathcal{F}$ represent a class of real-valued functions $f:\rchi\times\rchi\to\mathbb{R}$ uniformly bounded by $R$ with the envelope function $F$. Then there exists a universal constant $C > 0$ such that
\[
\mathbb{E}\left[\sup_{f\in\mathcal{F}}\left|\frac{\sum_{1\le i\neq j\le n}\epsilon_i\epsilon_jf_{i,j}(X_i, X_j)}{\sqrt{n(n-1)}}\right|\right] \le C\norm{F}_{2,P}J_2(a, \mathcal{F}, \norm{\cdot}_{2})\left[1+\frac{J_2(a, \mathcal{F}, \norm{\cdot}_2)b^2}{a^2}\right],
\]
for any $a \ge A_n$ and $b\ge B_n$, where $B_n^2 = R/(n\norm{F}_{2,P})$,
\begin{align*}
\norm{F}_{2,P}^2 &:= \frac{1}{n(n-1)}\sum_{1\le i\neq j\le n}\mathbb{E}[F^2_{i,j}(X_i, X_j)],\\
A_n^2 &:= \norm{F}_{2,P}^{-2}\left[\Gamma_{n,1}^2(\mathcal{F}) + \Gamma_{n,2}^2(\mathcal{F}) + \Sigma_n^2(\mathcal{F})\right],\\
\Gamma_{n,1}^2(\mathcal{F}) &:= \mathbb{E}\left[\sup_{f\in\mathcal{F}}\frac{1}{n(n-1)}\left|\sum_{1\le i\neq j\le n} \left\{\mathbb{E}\left[f^2_{i,j}(X_i, X_j)|X_i\right] - \mathbb{E}\left[f^2_{i,j}(X_i, X_j)\right]\right\}\right|\right],\\
\Gamma_{n,2}^2(\mathcal{F}) &:= \mathbb{E}\left[\sup_{f\in\mathcal{F}}\frac{1}{n(n-1)}\left|\sum_{1\le i\neq j\le n} \left\{\mathbb{E}\left[f^2_{i,j}(X_i, X_j)|X_j\right] - \mathbb{E}\left[f^2_{i,j}(X_i, X_j)\right]\right\}\right|\right],\\
\Sigma_{n}^2(\mathcal{F}) &:= \sup_{f\in\mathcal{F}}\frac{1}{n(n-1)}\sum_{1\le i\neq j\le n}\mathbb{E}\left[f^2_{i,j}(X_i, X_j)\right].
\end{align*}
\end{thm}

\begin{proof}
{\cred See Appendix~\ref{appsec:MaximalInequ} for a proof.}
\end{proof}

% \section{Examples Related to Non-parametric Statistics}\label{sec:Nonpar}
% \section{Concluding Remarks}\label{sec:Conclusion}

\newpage
\appendix
\setcounter{section}{0}
\setcounter{equation}{0}
\setcounter{figure}{0}
\renewcommand{\thesection}{S.\arabic{section}}
\renewcommand{\theequation}{E.\arabic{equation}}
\renewcommand{\thefigure}{A.\arabic{figure}}
% \tableofcontents
% \titlelabel{\thetitle: }
% \cftsetindents{section}{1em}{2.5em}
% \cftsetindents{subsection}{1.5em}{3em}
% \setcounter{page}{1}
  \begin{center}
  \Large {\bf Appendix to ``Tail bounds for canonical $U$-statistics and $U$-processes with unbounded kernels''}
  \end{center}

%%%%%%%%%%%%%%%%%%%%%%%%%%%%%%%%%%%%%%%%%%%%%%%%%%%%
%%%%%%%%%%%%%%%%%%%%%%%%%%%%%%%%%%%%%%%%%%%%%%%%%%%%
\section{Proofs of Results in Section~\ref{sec:UStat}}\label{sec:AppUStat}
\subsection{Proof of Theorem~\ref{thm:MainUStat-bounded-kernel}}\label{appsec:proof-thm:MainUStat-bounded-kernel}
By Theorem 3.5.3 of \cite{DeLaPena99}, it follows that
\[
\|U_n^D\|_p \le \mathfrak{C}\left\|\mathcal{U}_n\right\|_p,
\]
where
\[
\mathcal{U}_n := \sum_{1\le i\neq j\le n} f_{i,j}(Z_i, Z_j').
\]
For $1\le i\le n$ and any $z$, define
\[
h_i(z) := \sum_{\substack{1\le j\le n,\\j\neq i}}\, f_{i,j}(z, Z_j').
\]
Observe that
\[
\mathcal{U}_n = \sum_{i=1}^n h_i(Z_i).
\]
First, we consider the behavior of $h_i(z)$ for a fixed $z$ and then the behavior of $\mathcal{U}_n$. By the degeneracy of the kernel, we have that $h_i(z)$ is a sum of independent mean zero random variables for every $i, z$. Moreover, $\|f_{i,j}(z, Z_j')\|_{\psi_{\beta}} \le F_i(z)K_G$ for all $i, z$. Hence, Theorem 3.4 of~\cite{kuchibhotla2022moving} (with $q = 1$ and $t = \log(\delta/3)$) implies
\[
\mathbb{P}\left(|h_i(z)| \ge 7\sqrt{\log(3/\delta)}\left(\sum_{\substack{1\le j\le n,\\j\neq i}} \mathbb{E}[f_{i,j}^2(z, Z_j')]\right)^{1/2} + C_{\beta}F_i(z)K_G(\log(2n))^{1/\beta}(\log(3/\delta))^{1/\beta^*}\right) \le \delta,
\]
where $\beta^* = \min\{\beta, 1\}$. Based on this, define
\[
H_i(z; \delta) := F_i(z)\mathcal{B}_i(z, \delta/n),
\]
where
\[
\mathcal{B}_i(z, \delta) := 7\sqrt{\log(3/\delta)}\left(\sum_{\substack{1\le j\le n,\\j\neq i}} \mathbb{E}\left[\frac{f_{i,j}^2(z, Z_j')}{F_i^2(z)}\right]\right)^{1/2} + C_{\beta}K_G(\log(2n))^{1/\beta}(\log(3/\delta))^{1/\beta^*}.
\]

Getting back to the behavior of $\mathcal{U}_n$, we first note that by degeneracy and symmetrization,
\begin{equation}\label{eq:symmetrization}
\|\mathcal{U}_n\|_p \le 2\left\|\sum_{i=1}^n \varepsilon_ih_i(Z_i)\right\|_p\quad\mbox{for all}\quad p\ge1.
\end{equation}
Here $\varepsilon_1, \ldots, \varepsilon_n$ are independent Rademacher random variables independent of $Z_1, \ldots, Z_n$, $Z_1', \ldots, Z_n'$. Hence, it suffices to understand the behavior of
\[
\mathcal{U}_n' := \sum_{i=1}^n \varepsilon_ih_i(Z_i).
\]
(The introduction of Rademacher variables is only done for notational convenience in applying truncation.)
\begin{equation}\label{eq:main-part-1}
\begin{split}
    \mathbb{P}\left(|\mathcal{U}_n'| \ge t\right) &\le \mathbb{P}\left(\left|\sum_{i=1}^n \varepsilon_ih_i(Z_i)\mathbf{1}\{|h_i(Z_i)| \le H_i(Z_i; \delta_1)\}\right| \ge t\right)\\
&\quad+ \mathbb{P}\left(|h_i(Z_i)| > H_i(Z_i; \delta_1)\quad\mbox{for some}\quad 1\le i\le n\right)\\
    &\le \mathbb{P}\left(\left|\sum_{i=1}^n \varepsilon_ih_i(Z_i)\mathbf{1}\{|h_i(Z_i)| \le H_i(Z_i; \delta_1)\}\right| \ge t\right)\\
&\quad+ \sum_{i=1}^n \mathbb{P}(|h_i(Z_i)| > H_i(Z_i; \delta_1))\\
    &\le \mathbb{P}\left(\left|\sum_{i=1}^n \varepsilon_ih_i(Z_i)\mathbf{1}\{|h_i(Z_i)| \le H_i(Z_i; \delta_1)\}\right| \ge t\right) + \delta_1.
\end{split}
\end{equation}
Because $\{h_i(Z_i): 1\le i\le n\}$ are independent random variables conditional on $\{Z_j': 1\le j\le n\}$, we get by another application of Theorem 3.4 of~\cite{kuchibhotla2022moving} (with $q = 1$ and $t = \log(3/\delta_2)$) that conditional on $\{Z_j': 1\le i\le n\}$, with probability at least $1 - \delta_2$,
\begin{equation}\label{eq:main-part-2}
\begin{split}
&\left|\sum_{i=1}^n \varepsilon_ih_i(Z_i)\mathbf{1}\{|h_i(Z_i)| \le H_i(Z_i)\}\right|\\ &\quad\le 7\sqrt{\log(3/\delta_2)}\left(\sum_{i=1}^n \mathbb{E}[h_i^2(Z_i)|\{Z_j'\}]\right)^{1/2}\\
&\qquad+ C_{\alpha}(\log(2n))^{1/\alpha}(\log(3/\delta_2))^{1/\alpha^*}\max_{1\le i\le n}\left\|h_i(Z_i)\mathbf{1}\{|h_i(Z_i)| \le H_i(Z_i; \delta_1)\}\right\|_{\psi_{\alpha}|\{Z_j'\}}.
\end{split}
\end{equation}
Observe now that
\begin{align*}
    &\left\|h_i(Z_i)\mathbf{1}\{|h_i(Z_i)| \le H_i(Z_i; \delta_1)\}\right\|_{\psi_{\alpha}|\{Z_j'\}}\\ &\le \|F_i(Z_i)\mathcal{B}_i(Z_i, \delta_1/n)\|_{\psi_{\alpha}}\\
    &\le 7\sqrt{\log(3n/\delta_1)}\left\|\left(\sum_{\substack{1\le j\le n,\\j\neq i}} \mathbb{E}[f_{i,j}^2(Z_i, Z_j')|Z_i]\right)^{1/2}\right\|_{\psi_{\alpha}} + C_{\beta}K_FK_G(\log(2n))^{1/\beta}(\log(3n/\delta_1))^{1/\beta^*}.
\end{align*}
To bound the first term on the right hand side of~\eqref{eq:main-part-2}, we follow the argument in the proof of Theorem 3.2 of~\cite{Gine00}. However, in place of inequality (3.8) of~\cite{Gine00}, we apply Theorem B.1 of~\cite{kuchibhotla2022moving}. Following the display after inequality (3.11) of~\cite{Gine00}, we have
\begin{align*}
    &\left(\sum_{i=1}^n \mathbb{E}[h_i^2(Z_i)|\{Z_j'\}]\right)^{1/2}\\
    &= \sup\left\{\sum_{i=1}^n \mathbb{E}[h_i(Z_i)\gamma_i(Z_i)|\{Z_j'\}]:\,\sum_{i=1}^n \mathbb{E}[\gamma_i^2(Z_i)] \le 1\right\}\\
    &= \sup\left\{\sum_{j=1}^n \left(\sum_{\substack{1\le i\le n,\\i\neq j}} \mathbb{E}[f_{i,j}(Z_i, Z_j')\gamma_i(Z_i)|Z_j']\right):\, \sum_{i=1}^n \mathbb{E}[\gamma_i^2(Z_i)] \le 1\right\},
\end{align*}
where the supremum is taken over a countable subset of mean zero vector functions $(\gamma_1, \ldots, \gamma_n)$. Define
\[
W_j(\gamma) = \sum_{\substack{1\le i\le n,\\i\neq j}} \mathbb{E}[f_{i,j}(Z_i, Z_j')\gamma_i(Z_i)|Z_j'].
\]
Degeneracy of $\{f_{i,j}\}$ implies that $W_j$'s are mean zero independent random variables. Hence, by Theorem B.1 of~\cite{kuchibhotla2022moving}, we get
\begin{align*}
&\left(\mathbb{E}\left[\sup_{\gamma}\left|\sum_{j=1}^n W_j\right|^p\right]\right)^{1/p}\\
&\le 2\mathbb{E}\left[\sup_{\gamma}\left|\sum_{j=1}^n W_j\right|\right] + \sqrt{2p}\left(\sup_{\gamma}\sum_{j=1}^n \mathbb{E}[W_j^2]\right)^{1/2}+ C_{\beta}p^{1/\beta^*}\left\|\max_{1\le i\le n}\sup_{\gamma}W_j\right\|_{\psi_{\beta}}.
\end{align*}
Following the argument in Theorem 3.2 of~\cite{Gine00}, we get
\begin{align*}
    \mathbb{E}\left[\sup_{\gamma}\left|\sum_{j=1}^n W_j\right|\right] ~&\le~ \left(\sum_{1\le i\neq j\le n}\mathbb{E}[f_{i,j}^2(Z_i, Z_j')]\right)^{1/2},\\
    \sup_{\gamma}\sum_{j=1}^n \mathbb{E}[W_j^2] ~&\le~ \|(f_{i,j})\|_{L^2\to L^2}^2,\\
\max_{1\le j\le n}\sup_{\gamma}\,W_j ~&\le~ \max_{1\le j\le n}\left(\sum_{\substack{1\le i\le n,\\i\neq j}} \mathbb{E}[f_{i,j}^2(Z_i, Z_j')|Z_j']\right)^{1/2}.
\end{align*}
Therefore, by Markov's inequality, with probability at least $1 - \delta_3$,
\begin{equation}\label{eq:main-part-3}
\begin{split}
&\left(\sum_{i=1}^n \mathbb{E}[h_i^2(Z_i)|\{Z_j'\}]\right)^{1/2}\\ &\le 2\left(\sum_{1\le i\neq j\le n}\mathbb{E}[f_{i,j}^2(Z_i, Z_j')]\right)^{1/2} + \sqrt{2\log(1/\delta_3)}\|(f_{i,j})\|_{L^2\to L^2}\\
&\quad+ C_{\beta}(\log(1/\delta_3))^{1/\beta^*}(\log(n))^{1/\beta}\max_{1\le j\le n}\left\|\left(\sum_{\substack{1\le i\le n,\\i\neq j}}\mathbb{E}[f_{i,j}^2(Z_i, Z_j')|Z_j']\right)^{1/2}\right\|_{\psi_{\beta}}.
\end{split}
\end{equation}

Combining inequalities~\eqref{eq:main-part-1},~\eqref{eq:main-part-2},~\eqref{eq:main-part-3}, we get that with probability $1 - \delta_1 - \delta_2 - \delta_3$,
\begin{align*}
    |\mathcal{U}_n'| &\le 14\sqrt{\log(3/\delta_2)}\left(\sum_{i\neq j}\mathbb{E}[f_{i,j}^2(Z_i, Z_j')]\right)^{1/2}\\
    &\quad+ 7\sqrt{2\log(3/\delta_2)\log(1/\delta_3)}\|(f_{i,j})\|_{L^2\to L^2}\\
    &\quad+ \mathfrak{C}_{\beta}(\log(3/\delta_2))^{1/2}(\log(1/\delta_3))^{1/\beta^*}(\log n)^{1/\beta}\max_{1\le j\le n}\left\|\left(\sum_{\substack{1\le i\le n,\\i\neq j}}\mathbb{E}[f_{i,j}^2(Z_i, Z_j')|Z_j']\right)^{1/2}\right\|_{\psi_{\beta}}\\
    &\quad+ \mathfrak{C}_{\alpha}(\log(3n/\delta_1))^{1/2}(\log(3/\delta_2))^{1/\alpha^*}(\log(2n))^{1/\alpha}\max_{1\le i\le n}\left\|\left(\sum_{\substack{1\le j\le n,\\j\neq i}} \mathbb{E}[f_{i,j}^2(Z_i, Z_j')|Z_i]\right)^{1/2}\right\|_{\psi_{\alpha}}\\
    &\quad+ \mathfrak{C}_{\alpha,\beta}K_FK_G(\log(2n))^{1/\alpha + 1/\beta}(\log(3/\delta_2))^{1/\alpha^*}(\log(3n/\delta_1))^{1/\beta^*}.
\end{align*}
Taking $\delta_1 = \delta_2 = \delta_3 = \delta/3$ and integrating over $\delta\in[0, 1]$, this inequality yields the following moment bound
\begin{align*}
    \|\mathcal{U}_n'\|_p &\le \mathfrak{C}p^{1/2}\left(\sum_{i\neq j}\mathbb{E}[f_{i,j}^2(Z_i, Z_j')]\right)^{1/2}\\
    &\quad+ \mathfrak{C}p\|(f_{i,j})\|_{L^2\to L^2}\\
    &\quad+ \mathfrak{C}_{\beta}p^{1/2 + 1/\beta^*}(\log n)^{1/\beta}\max_{1\le j\le n}\left\|\left(\sum_{\substack{1\le i\le n,\\i\neq j}}\mathbb{E}[f_{i,j}^2(Z_i, Z_j')|Z_j']\right)^{1/2}\right\|_{\psi_{\beta}}\\
    &\quad+ \mathfrak{C}_{\alpha}p^{1/2 + 1/\alpha^*}(\log(2n))^{1/2 + 1/\alpha}\max_{1\le i\le n}\left\|\left(\sum_{\substack{1\le j\le n,\\j\neq i}} \mathbb{E}[f_{i,j}^2(Z_i, Z_j')|Z_i]\right)^{1/2}\right\|_{\psi_{\alpha}}\\
    &\quad+ \mathfrak{C}_{\alpha,\beta}p^{1/\alpha^* + 1/\beta^*}K_FK_G(\log(2n))^{1/\alpha + 1/\beta + 1/\beta^*}.
\end{align*}
This inequality combined with~\eqref{eq:symmetrization} yields the tail bound for $U_n^D$.

\subsection{Proof of Lemma~\ref{lem:ReductionSymmetrization}}\label{appsec:proof-lem:ReductionSymmetrization}
From Theorem 3.1.1 of \cite{DeLaPena99} and following the arguments similar to those in Theorem 3.5.3 of \cite{DeLaPena99}, we get for all $p \ge 1$
\[
\norm{T_n^D}_p \le 48\norm{\sum_{1\le i\neq j\le n}\varepsilon_i\varepsilon_j'f_{i,j}^D(Z_i, Z_j')}_p,
\]
where $\varepsilon_i, \varepsilon_i', 1\le i\le n$ are Rademacher random variables independent of $(Z_i, Z_i'), 1\le i\le n$. Note from \eqref{eq:SplitDegenerate} that
\begin{align*}
\varepsilon_i\varepsilon_j'f_{i,j}^D(Z_i, Z_j') &= \varepsilon_i\varepsilon_j'f_{i,j}(Z_i, Z_j') - \varepsilon_i\varepsilon_j'\int f_{i,j}(z, Z_j')P_i(dz)\\
&\quad- \varepsilon_i\varepsilon_j'\int f_{i,j}(Z_i, z)P_j(dz) + \varepsilon_i\varepsilon_j'\iint f_{i,j}(z, z')P_i(dz)P_j(dz).
\end{align*}
Here $P_i$ represents the probability measure of $Z_i$ for $1\le i\le n$. By Jensen's inequality, it is clear that for $p\ge 1$,
\begin{align*}
\norm{\sum_{1\le i\neq j\le n}\varepsilon_i\varepsilon_j'\int f_{i,j}(z, Z_j')P_i(dz)}_p &\le \norm{\sum_{1\le i\neq j\le n}\varepsilon_i\varepsilon_j'f_{i,j}(Z_i, Z_j')}_p,\\
\norm{\sum_{1\le i\neq j\le n}\varepsilon_i\varepsilon_j'\int f_{i,j}(Z_i, z)dP_j(z)}_p&\le \norm{\sum_{1\le i\neq j\le n}\varepsilon_i\varepsilon_j'f_{i,j}(Z_i, Z_j')}_p,\\
\norm{\sum_{1\le i\neq j\le n}\varepsilon_i\varepsilon_j'\iint f_{i,j}(z, z')P_i(dz)P_j(dz)}_p &\le \norm{\sum_{1\le i\neq j\le n}\varepsilon_i\varepsilon_j'f_{i,j}(Z_i, Z_j')}_p.
\end{align*}
Therefore, for $p\ge 1$,
\[
\norm{T_n^D}_p \le 192\norm{\sum_{1\le i\neq j\le n}\varepsilon_i\varepsilon_j'f_{i,j}(Z_i, Z_j')}_p
\]

Throughout the proofs in all the appendices to follow, we use the notation
\[
\mathcal{Z}_n' := \{(Z_1', \varepsilon_1'), \ldots, (Z_n', \varepsilon_n')\}\quad\mbox{and}\quad \mathcal{Z}_n := \{(Z_1, \varepsilon_1), \ldots, (Z_n, \varepsilon_n)\}.
\]
Note that this is different from $\mathcal{Z}_n'$ and $\mathcal{Z}_n$ defined in the main text.
\subsection{Proof of Theorem~\ref{thm:MainUStat}}\label{subsec:AppSecMainUStat}
Based on the basic decomposition~\eqref{eq:BasicDecomposition}, we get
\[
\sum_{1\le i\neq j\le n}\varepsilon_i\phi(Z_i)w_{i,j}(X_i, X_j')\psi(Z_j')\varepsilon_j' = \mathcal{U}_{n}^{(1)} + \mathcal{U}_n^{(2)} + \mathcal{U}_n^{(3)} + \mathcal{U}_n^{(4)},
\]
where
\begin{equation}\label{eq:SplitDegenerateUStats}
\begin{split}
\mathcal{U}_n^{(1)} &:= \sum_{1\le i\neq j\le n}\varepsilon_i\Phi_{i,1}w_{i,j}(X_i, X_j')\Psi'_{j,1}\varepsilon_j',\\
\mathcal{U}_n^{(2)} &:= \sum_{1\le i\neq j\le n}\varepsilon_i\Phi_{i,2}w_{i,j}(X_i, X_j')\Psi'_{j,1}\varepsilon_j',\\
\mathcal{U}_n^{(3)} &:= \sum_{1\le i\neq j\le n}\varepsilon_i\Phi_{i,1}w_{i,j}(X_i, X_j')\Psi'_{j,2}\varepsilon_j',\\
\mathcal{U}_n^{(4)} &:= \sum_{1\le i\neq j\le n}\varepsilon_i\Phi_{i,2}w_{i,j}(X_i, X_j')\Psi'_{j,2}\varepsilon_j'.
\end{split}
\end{equation}
It is easy to verify that $\mathcal{U}_n^{(k)}, 1\le k\le 4$ are all degenerate $U$-statistics. From Theorem 3.2 of \cite{Gine00}, we get that there exists a constant $K > 0$ such that for all $p\ge 1$,
\[
\norm{\mathcal{U}_{n}^{(1)}}_p \le K\left[\sqrt{p}A + pB + p^{3/2}C + p^{2}D\right],
\]
where
\begin{equation}\label{eq:DefinitionsABCD}
\begin{split}
A &:= \left(\sum_{1\le i\neq j\le n} \mathbb{E}\left[\Phi_{i,1}^2w_{i,j}^2(X_i, X_j')\left(\Psi_{i,1}'\right)^2\right]\right)^{1/2},\\
B &:= \sup\left\{\mathbb{E}\sum_{1\le i\neq j\le n} \varepsilon_i\xi_i(\varepsilon_i, Z_i)\Phi_{i,1}w_{i,j}(X_i, X_j')\Psi'_{i,1}\zeta_j(\varepsilon_j', Z_j')\varepsilon_j':\right.\\&\qquad\qquad\quad\left.\mathbb{E}\sum_{i = 1}^n \xi_i^2(\varepsilon_i, Z_i) \le 1,\, \mathbb{E}\sum_{j = 1}^n \zeta_i^2(\varepsilon_i', Z_i') \le 1\right\},\\
C^p &:= \mathbb{E}\left(\max_{1\le i\le n}\mathbb{E}\left[\sum_{j = 1}^n \Phi_{i,1}^2w_{i,j}^2(X_i, X_j')\left(\Psi'_{i,1}\right)^2\big|X_i, Y_i\right]\right)^{p/2}\\ &\qquad+ \mathbb{E}\left(\max_{1\le j\le n}\mathbb{E}\left[\sum_{i = 1}^n \Phi_{i,1}^2w_{i,j}^2(X_i, X_j')\left(\Psi'_{i,1}\right)^2\big|X_j', Y_j'\right]\right)^{p/2}\\
D^p &:= \mathbb{E}\left[\max_{1\le i\neq j \le n} |\Phi_{i,1}w_{i,j}(X_i, X_j')\Psi_{j,1}'|^{p}\right].
\end{split}
\end{equation}
It is clear that
\begin{align*}
A^2 \le \sum_{1\le i\neq j\le n}\mathbb{E}\left[\phi^2(Y_i)w_{i,j}^2(X_i, X_j)\psi^2(Y_j)\right]= \sum_{1\le i\neq j\le n}\mathbb{E}\left[\sigma_{i,\phi}^2(X_i)w_{i,j}^2(X_i, X_j)\sigma_{j,\psi}^2(X_j)\right].
\end{align*}
The quantity $B$ appears as the square root of the wimpy variance of the supremum of an empirical process; see \citet[page 314]{MR3185193}. Lemma~\ref{lem:BoundingB} of Section~\ref{sec:AuxLemmas} implies that
\begin{align*}
B &\le\sup\left\{\sum_{1\le i\neq j\le n}\mathbb{E}\left[q_i(X_i)\sigma_{i,\phi}(X_i)w_{i,j}(X_i, X_j)\sigma_{j,\psi}(X_j)p_j(X_j)\right]:\right.\\
&\qquad\qquad\quad\left.\sum_{j = 1}^n \mathbb{E}\left[q_i^2(X_i)\right] \le 1,\,\sum_{i=1}^n \mathbb{E}\left[p_j^2(X_j)\right] \le 1\right\}.
\end{align*}
For bounding $C$, note that
\begin{align*}
\mathbb{E}\left[\sum_{j = 1}^n \Phi_{i,1}^2w_{i,j}^2(X_i, X_j')\left(\Psi'_{j,1}\right)^2\big|X_i, Y_i\right] &\le T_{\phi}^2\sup_{x}\sum_{j = 1}^n \mathbb{E}\left[w_{i,j}^2(x, X_j)\sigma^2_{j,\psi}(X_j)\right],\\
\mathbb{E}\left[\sum_{i = 1}^n \Phi_{i,1}^2w_{i,j}^2(X_i, X_j')\left(\Psi'_{j,1}\right)^2\big|X_j', Y_j'\right] &\le T_{\psi}^2\sup_{x}\sum_{i = 1}^n \mathbb{E}\left[w_{i,j}^2(X_i, x)\sigma^2_{i, \phi}(X_i)\right].
\end{align*}
Combining these two inequalities implies that
\[
C \le T_{\phi}\sup_{x}\left(\sum_{j = 1}^n \mathbb{E}\left[w_{i,j}^2(x, X_j)\sigma^2_{j,\psi}(X_j)\right]\right)^{1/2} + T_{\psi}\sup_{x}\left(\sum_{i = 1}^n \mathbb{E}\left[w_{i,j}^2(X_i, x)\sigma^2_{i, \phi}(X_i)\right]\right)^{1/2}
\]
Finally, it is clear from assumption~\ref{assump:BoundedKernels} that $D \le T_{\phi}T_{\psi}B_{w}.$ Combining all these with Theorem 3.2 of~\cite{Gine00} and noting
\[
T_{\phi} \le K_{\alpha}C_{\phi}(\log n)^{1/\alpha}\quad\mbox{and}\quad T_{\psi} \le K_{\beta} C_{\psi}(\log n)^{1/\beta},
\]
we get that there exists a constant $K > 0$ such that for all $p\ge 1$
\begin{equation}\label{eq:Un1Bound}
\norm{\mathcal{U}_n^{(1)}}_p \le K\left[\sqrt{p}\Lambda_{1/2} + p\Lambda_1 + p^{3/2}\left\{\Lambda_{3/2}^{(\alpha)} + \Lambda_{3/2}^{(\beta)}\right\} + p^2\Lambda_{2}\right].
\end{equation}
To bound $\mathcal{U}_n^{(2)}$ and $\mathcal{U}_n^{(3)}$ in~\eqref{eq:SplitDegenerateUStats}, we use Hoffmann-J{\o}gensen's inequality (Proposition 6.8 of \cite{LED91}). Observe that
\[
\mathcal{U}_n^{(2)} := \sum_{i = 1}^n \varepsilon_i\Phi_{i,2}g_i(X_i; \mathcal{Z}_n'),\quad\mbox{where}\quad g_i(X_i; \mathcal{Z}_n') := \sum_{j = 1, j\neq i}^n w_{i,j}(X_i, X_j')\Psi_{j,1}'\varepsilon_j'.
\]
With $\mathcal{Z}_n' := \{(\varepsilon_1', Z_1'),\ldots,(\varepsilon_n', Z_n')\}$ and $\mathcal{X}_n := \{X_1, \ldots, X_n\}$, note that
\begin{align*}
\mathbb{P}&\left(\max_{1\le I\le n}\left|\sum_{i = 1}^I \varepsilon_i\Phi_{i,2}g_i(X_i, \mathcal{Z}_n')\right| > 0\big|\mathcal{X}_n, \mathcal{Z}_n'\right)%\\ &
\le \mathbb{P}\left(\max_{1\le i\le n}|\phi(Z_i)| \ge T_{\phi}\big|\mathcal{X}_n\right) \le 1/8,
\end{align*}
and so, by Equation (6.8) of~\cite{LED91}, we get
\begin{align*}
\mathbb{E}\left[\mathcal{U}_n^{(2)}\big|\mathcal{X}_n, \mathcal{Z}_n'\right] &\le 8\mathbb{E}\left[\max_{1\le i\le n}\left|\Phi_{i,2}\left(g_i(X_i; \mathcal{Z}_n')\right)\right|\big|\mathcal{X}_n, \mathcal{Z}_n'\right]\\
&\le 8\mathbb{E}\left[\max_{1\le i\le n}\left|\phi(Z_i)\right|\big|\mathcal{X}_n\right]\max_{1\le i\le n}\left|g_i(X_i; \mathcal{Z}_n')\right| = T_{\phi}\max_{1\le i\le n}\left|g_i(X_i; \mathcal{Z}_n')\right|.
\end{align*}
From assumption~\ref{assump:ConditionalTailAssumption} and Theorem 6.21 of~\cite{LED91}, we thus get for $0 < \alpha \le 1$,
\begin{equation}\label{eq:LTApplications}
\begin{split}
\norm{\mathcal{U}_n^{(2)}}_{\psi_{\alpha}\big|\mathcal{X}_n, \mathcal{Z}_n'} &\le K_{\alpha}\mathbb{E}\left[\mathcal{U}_n^{(2)}\big|\mathcal{X}_n, \mathcal{Z}_n'\right] + K_{\alpha}\norm{\max_{1\le i\le n}\left|\Phi_{i,2}g_i(X_i; \mathcal{Z}_n')\right|}_{\psi_{\alpha}\big|\mathcal{X}_n, \mathcal{Z}_n'}\\
&\le K_{\alpha}\left(T_{\phi} + \norm{\max_{1\le i\le n}|\phi(Y_i)|}_{\psi_{\alpha}\big|\mathcal{X}_n, \mathcal{Z}_n'}\right)\max_{1\le i\le n}\left|g_i(X_i; \mathcal{Z}_n')\right|\\
&\le K_{\alpha}C_{\phi}(\log n)^{1/\alpha}\max_{1\le i\le n}\left|g_i(X_i; \mathcal{Z}_n')\right|,
\end{split}
\end{equation}
for some constant $K_{\alpha}$ depending only on $\alpha$ (and can be different in different lines). If $\alpha \ge 1$, then we get
\[
\norm{\mathcal{U}_n^{(2)}}_{\psi_{\alpha^*}\big|\mathcal{X}_n, \mathcal{Z}_n'} \le K_{\alpha}C_{\phi}(\log n)^{1/\alpha}\max_{1\le i\le n}\left|g_i(X_i; \mathcal{Z}_n')\right|.
\]
See proof of Theorem 3.3 in~\cite{kuchibhotla2022moving} for similar argument. Thus,
\[
\mathbb{E}\left[|\mathcal{U}_n^{(2)}|^p\big|\mathcal{X}_n, \mathcal{Z}_n'\right] \le K_{\alpha}^pC_{\phi}^p(\log n)^{p/\alpha}p^{p/\alpha^*}\max_{1\le i\le n}\left|g_i(X_i; \mathcal{Z}_n')\right|^p.
\]
Thus, for $p\ge 1$,
\begin{equation}\label{eq:Un2Bound}
\mathbb{E}\left[|\mathcal{U}_n^{(2)}|^p\right] \le K_{\alpha}^pC_{\phi}^p(\log n)^{p/\alpha}p^{p/\alpha^*}\mathbb{E}\left[\max_{1\le i\le n}\left|g_i(X_i; \mathcal{Z}_n')\right|^p\right].
\end{equation}
To control the right hand side above, recall that
\[
g_i(x; \mathcal{Z}_n') = \sum_{j = 1, j\neq i}^n w_{i,j}(x, X_j')\Psi_{j,1}'\varepsilon_j',
\]
is a sum of mean zero independent random variables that are bounded by $B_wT_{\psi}$. Also, note that
\[
\mbox{Var}(g_i(x; \mathcal{Z}_n')) = \sum_{j=1,j\neq i}^n \mathbb{E}\left[w_{i,j}^2(x, X_j')\psi^2(Z_j')\right] = \sum_{j = 1, j\neq i}^n \mathbb{E}[w_{i,j}^2(x, X_j')\sigma_{j,\psi}^2(X_j')].
\]
Therefore by Bernstein's inequality (Lemma 4 of \cite{Geer13}), we get that
\begin{equation}\label{eq:UnionBernstein}
\mathbb{P}\left(\max_{1\le i\le n}|g_i(X_i; \mathcal{Z}_n')| - \Upsilon_{\psi}\sqrt{6\log(1 + n)} - 3B_wT_{\psi}\log n \ge \Upsilon_{\psi}\sqrt{t} + 3B_wT_{\psi}t\right) \le 2e^{-t},
\end{equation}
where
\[
\Upsilon_{\psi}^2 := \max_{x}\sum_{j = 1, j\neq i}^n \mathbb{E}[w_{i,j}^2(x, X_j')\sigma_{j,\psi}^2(X_j')].
\]
So, by Propositions A.3 and A.4 of~\cite{kuchibhotla2022moving}, we get that for $p\ge 1$,
\[
\mathbb{E}\left[\max_{1\le i\le n}|g_i(X_i; \mathcal{Z}_n')|^p\right] \le C^p\left[(\log n)^{p/2}\Upsilon_{\psi}^{p} + (B_wT_{\psi})^p(\log n)^p + p^{p/2}\Upsilon_{\psi}^p + p^p(B_wT_{\psi})^p\right].
\]
Hence for $p\ge 1$,
\begin{equation}\label{eq:Un2Bound}
\begin{split}
\mathbb{E}\left[|\mathcal{U}_n^{(2)}|^p\right] &\le K_{\alpha}^pC_{\phi}^p(\log n)^{p/\alpha}p^{p/\alpha^*}\left[(\log n)^{p/2}\Upsilon_{\psi}^{p} + (B_wT_{\psi})^p(\log n)^p\right]\\
&\quad+ K_{\alpha}^pC_{\phi}^p(\log n)^{p/\alpha}p^{p/\alpha^*}\left[p^{p/2}\Upsilon_{\psi}^p + p^p(B_wT_{\psi})^p\right].
\end{split}
\end{equation}
A similar calculation for $\mathcal{U}_n^{(3)}$ shows that for $p\ge 1$,
\begin{equation}\label{eq:Un3Bound}
\begin{split}
\mathbb{E}\left[|\mathcal{U}_n^{(3)}|^p\right] &\le K_{\beta}^pC_{\psi}^p(\log n)^{p/\beta}p^{p/\beta^*}\left[(\log n)^{p/2}\Upsilon_{\phi}^{p} + (B_wT_{\phi})^p(\log n)^p\right]\\
&\quad+ K_{\beta}^pC_{\psi}^p(\log n)^{p/\beta}p^{p/\beta^*}\left[p^{p/2}\Upsilon_{\phi}^p + p^p(B_wT_{\phi})^p\right],
\end{split}
\end{equation}
where
\[
\Upsilon_{\phi}^2 := \max_{x}\sum_{i = 1, i\neq j}^n \mathbb{E}[w_{i,j}^2(X_i, x)\sigma_{i,\phi}^2(X_i)].
\]
To control $\mathcal{U}_n^{(4)}$, recall that
\[
\mathcal{U}_n^{(4)} = \sum_{i=1}^n \varepsilon_i\Phi_{i,2}\left(\sum_{j = 1, j\neq i}^n w_{i,j}(X_i, X_j')\Psi_{j,2}'\varepsilon_j'\right).
\]
Following the arguments leading to~\eqref{eq:LTApplications}, we have
\[
\norm{\mathcal{U}_n^{(4)}}_{\psi_{\alpha^*}|\mathcal{X}_{n}, \mathcal{Z}_n'} \le K_{\alpha}C_{\phi}(\log n)^{1/\alpha}\max_{1\le i\le n}\left|\sum_{j = 1, j\neq i}^n w_{i,j}(X_i, X_j')\Psi_{j}'\varepsilon_j'\right|.
\]
Conditioning on $\mathcal{X}_n, \mathcal{X}_n'$, the right hand side satisfies the hypothesis of (6.8) of~\cite{LED91} and so by Theorem 6.21 of~\cite{LED91}, we get
\begin{align*}
\norm{\max_{1\le i\le n}\left|\sum_{j = 1, j\neq i}^n w_{i,j}(X_i, X_j')\Psi_{j}'\varepsilon_j'\right|}_{\psi_{\beta^*}|\mathcal{X}_n, \mathcal{X}_n'} &\le K_{\beta}\norm{\max_{1\le i\neq j\le n}\left|w_{i,j}(X_i, X_j')\psi(Z_j')\right|}_{\psi_{\beta}|\mathcal{X}_n, \mathcal{X}_n'}\\
&\le K_{\beta}B_{w}(\log n)^{1/\beta}C_{\psi},
\end{align*}
for some constant $K_{\beta}$ depending only on $\beta$. Therefore, for $p\ge 1$
\begin{equation}\label{eq:Un4Bound}
\mathbb{E}\left[|\mathcal{U}_n^{(4)}|^p\right] \le K^pC_{\psi}^pC_{\phi}^pp^{p(1/\alpha^{*} + 1/\beta^{*})}(\log n)^{p(\alpha^{-1} + \beta^{-1})}B_w^p \le K^pp^{p(1/\alpha^{*} + 1/\beta^{*})}\Lambda_2^p,
\end{equation}
for some constant $K > 0$.

Combining bounds~\eqref{eq:Un2Bound} and~\eqref{eq:Un3Bound}, we get that for some constant $K > 0$ and for all $p\ge 1$,
\begin{align*}
\norm{\mathcal{U}_n^{(2)} + \mathcal{U}_n^{(3)}}_p &\le Kp^{1/\alpha^*}(\log n)^{1/2}\Lambda_{3/2}^{(\alpha)} + Kp^{1/\beta^*}(\log n)^{1/2}\Lambda_{3/2}^{(\beta)}\\ &\qquad+ Kp^{1/2 + 1/\alpha^*}\Lambda_{3/2}^{(\alpha)} + Kp^{1/2 + 1/\beta^*}\Lambda_{3/2}^{(\beta)}\\ &\qquad+ K(\log n)\Lambda_2[p^{1/\alpha^*} + p^{1/\beta^*}] + K\Lambda_2[p^{1 + 1/\alpha^*} + p^{1 + 1/\beta^*}].
\end{align*}
Combining this inequality with~\eqref{eq:Un1Bound} and~\eqref{eq:Un4Bound}, we get for all $p\ge 1$
\begin{align*}
\norm{\sum_{\ell = 1}^4\mathcal{U}_n^{(\ell)}}_p &\le K\left[p^{1/2}\Lambda_{1/2} + p\Lambda_1 + p^{3/2}\left\{\Lambda_{3/2}^{(\alpha)} + \Lambda_{3/2}^{(\beta)}\right\} + p^2\Lambda_2\right]\\
&\quad+ Kp^{1/\alpha^*}(\log n)^{1/2}\Lambda_{3/2}^{(\alpha)} + Kp^{1/\beta^*}(\log n)^{1/2}\Lambda_{3/2}^{(\beta)}\\ &\quad+ Kp^{1/2 + 1/\alpha^*}\Lambda_{3/2}^{(\alpha)} + Kp^{1/2 + 1/\beta^*}\Lambda_{3/2}^{(\beta)}\\
&\quad+ K(\log n)\Lambda_2[p^{1/\alpha^*} + p^{1/\beta^*}] + K\Lambda_2[p^{1 + 1/\alpha^*} + p^{1 + 1/\beta^*}]\\
&\quad+ Kp^{(1/\alpha^{*} + 1/\beta^{*})}\Lambda_2.
\end{align*}
Since $\alpha^* \le 1$ and $\beta^* \le 1$, we have
$$\min\{p^{1/2 + 1/\alpha^*}, p^{1/2 + 1/\beta^*}\} \ge p^{3/2}\mbox{ and }\min\{p^{1 + 1/\alpha^*}, p^{1 + 1/\beta^*}, p^{1/\alpha^* + 1/\beta^*}\} \ge p^2.$$ Using these inequalities, the bound above can be simplified as
\begin{align*}
\norm{\sum_{\ell = 1}^4\mathcal{U}_n^{(\ell)}}_p &\le K p^{1/2}\Lambda_{1/2} + K p\Lambda_1\\
&\quad + K p^{1/\alpha^*}\left[(\log n)^{1/2}\Lambda_{3/2}^{(\alpha)} + (\log n)\Lambda_2\right] + Kp^{1/\beta^*}\left[(\log n)^{1/2}\Lambda_{3/2}^{(\beta)} + (\log n)\Lambda_2\right]\\
&\quad + K p^{1/2 + 1/\alpha^*}\Lambda_{3/2}^{(\alpha)} + K p^{1/2 + 1/\beta^*}\Lambda_{3/2}^{(\beta)}\\
&\quad + Kp^{1/\alpha^* + 1/\beta^*}\Lambda_2.
\end{align*}
Here the constant $K > 0$ depends only on $\alpha, \beta$. This completes the proof based on Lemma~\ref{lem:ReductionSymmetrization}.
%%%%%%%%%%%%%%%%%%%%%%%%%%%%%%%%%%%%%%%%%%%%%%%%%%
%%%%%%%%%%%%%%%%%%%%%%%%%%%%%%%%%%%%%%%%%%%%%%%%%%

\subsection{Auxiliary Lemmas Used in Theorem~\ref{thm:MainUStat}}\label{sec:AuxLemmas}

The two lemmas to follow in this section provide explicit (but not necessarily optimal) constants for Equations (3.1) and (2.6) of~\cite{Gine00}. These lemmas can be used in the proof of Theorem 3.2 of~\cite{Gine00} to get explicit constants. In this respect, we note that Theorem 3.4.8 of~\cite{GINE16} (which was first proved in~\cite{HOUDRE03}) does not imply Theorem 3.2 of~\cite{Gine00} since the result of~\cite{Gine00} applies for unbounded kernels in $U$-statistics while the result of~\cite{GINE16} applies exclusively for bounded kernel $U$-statistics.
\begin{lem}\label{lem:BernsteinGeneral}
Suppose $Z_1, \ldots, Z_n$ are independent mean zero random variables. Then for $p\ge 1$,
\[
\mathbb{E}\left[\left|\sum_{i=1}^n Z_i\right|^p\right] \le 4^pp^{p/2}\left(\sum_{i=1}^n \mathbb{E}\left[Z_i^2\right]\right)^{p/2} + 4^pp^p\mathbb{E}\left[\max_{1\le i\le n}|Z_i|^p\right].
\]
\end{lem}
\begin{proof}
By Theorem 7 of~\cite{Bouch05}, we get for $p\ge 2$,
\[
\mathbb{E}\left[\left|\sum_{i=1}^n Z_i\right|^p\right] \le 2^{p+1}\left(\frac{2p}{e - \sqrt{e}}\right)^{p/2}\mathbb{E}\left[\left(\sum_{i=1}^n Z_i^2\right)^{p/2}\right].
\]
By Theorem 8 of~\cite{Bouch05}, we get for $p\ge 2$,
\[
\mathbb{E}\left[\left(\sum_{i=1}^n Z_i^2\right)^{p/2}\right] \le 3^{p/2}\left(\sum_{i=1}^n \mathbb{E}\left[Z_i^2\right]\right)^{p/2} + \left(\frac{3p\kappa}{2}\right)^{p/2}\mathbb{E}\left[\max_{1\le i\le n}|Z_i|^{p}\right],
\]
for $\kappa = 0.5\sqrt{e}/(\sqrt{e} - 1)$. Thus for $p\ge 2$,
\[
\mathbb{E}\left[\left|\sum_{i=1}^n Z_i\right|^p\right] \le 4^pp^{p/2}\left(\sum_{i=1}^n \mathbb{E}\left[Z_i^2\right]\right)^{p/2} + 4^pp^p\mathbb{E}\left[\max_{1\le i\le n}|Z_i|^p\right].
\]
Since the inequality holds true for $p = 1$ trivially, the result follows.
\end{proof}
\begin{lem}\label{lem:Gine2p6}
Suppose $\xi_i, 1\le i\le n$ are independent random variables, then for $p\ge 1$ and $\alpha > 0$,
\[
p^{p\alpha}\sum_{i=1}^n \mathbb{E}\left[\left|\xi_i\right|^p\right] \le 4(1.5)^{p\alpha}p^{p\alpha}\mathbb{E}\left[\max_{1\le i\le n}|\xi_i|^p\right] + 2(1.5)^{p\alpha}\left(\sum_{i=1}^n \mathbb{E}\left[|\xi_i|\right]\right)^p.
\]
\end{lem}
\begin{proof}
Fix $p\ge 1$. Define $\delta_0 \ge 0$ such that
\[
\delta_0 := \inf\left\{t > 0:\,\sum_{i=1}^n \mathbb{P}\left(|\xi_i| > t\right) \le 1\right\}.
\]
By (1.4.4) of \cite{DeLaPena99}, it follows that
\begin{equation}\label{eq:144DeLaGine}
\frac{1}{2}\max\left\{\delta_0^p, \sum_{i=1}^n \mathbb{E}\left[|\xi_i|^p\mathbbm{1}_{\{|\xi_i| > \delta_0\}}\right]\right\} \le \mathbb{E}\left[\max_{1\le i\le n}|\xi_i|^p\right].
\end{equation}
Observe that
\begin{align*}
\sum_{i=1}^n \mathbb{E}\left[|\xi_i|^p\right] &= \sum_{i=1}^n \mathbb{E}\left[|\xi_i|^p\mathbbm{1}_{\{|\xi_i| > \delta_0\}}\right] + \sum_{i=1}^n \mathbb{E}\left[|\xi_i|^p\mathbbm{1}_{\{|\xi_i| \le \delta_0\}}\right]\\
&\overset{(a)}{\le} 2\mathbb{E}\left[\max_{1\le i\le n}|\xi_i|^p\right] + \sum_{i=1}^n \mathbb{E}\left[|\xi_i|^p\mathbbm{1}_{\{|\xi_i| \le \delta_0\}}\right]\\
&\le 2\mathbb{E}\left[\max_{1\le i\le n}|\xi_i|^p\right] + \delta_0^{p-1}\sum_{i=1}^n \mathbb{E}\left[|\xi_i|\mathbbm{1}_{\{|\xi_i| \le \delta_0\}}\right]\\
&\overset{(a)}{\le} 2\mathbb{E}\left[\max_{1\le i\le n}|\xi_i|^p\right] + 2\mathbb{E}\left[\max_{1\le i\le n}|\xi_i|^{p-1}\right]\left(\sum_{i=1}^n \mathbb{E}\left[|\xi_i|\mathbbm{1}_{\{|\xi_i| \le \delta_0\}}\right]\right)\\
&\le 2\mathbb{E}\left[\max_{1\le i\le n}|\xi_i|^p\right] + 2\mathbb{E}\left[\max_{1\le i\le n}|\xi_i|^{p-1}\right]\left(\sum_{i=1}^n \mathbb{E}\left[|\xi_i|\right]\right).
\end{align*}
Inequality (a) follows from~\eqref{eq:144DeLaGine}.
To prove the result now, we consider two cases:
\begin{itemize}
	\item[--] \textit{Case 1:} If $$p^{p\alpha}\mathbb{E}\left[\max_{1\le i\le n}|\xi_i|^p\right] \le \left(\sum_{i=1}^n \mathbb{E}[|\xi_i|]\right)^p,$$ then
	\begin{align*}
	\mathbb{E}\left[\max_{1\le i\le n}|\xi_i|^{p-1}\right]\left(\sum_{i=1}^n \mathbb{E}\left[|\xi_i|\right]\right) &\le \left(\mathbb{E}\left[\max_{1\le i\le n}|\xi_i|^p\right]\right)^{(p-1)/p}\left(\sum_{i=1}^n \mathbb{E}\left[|\xi_i|\right]\right)\\
	&\le \frac{1}{p^{(p-1)\alpha}}\left(\sum_{i=1}^n \mathbb{E}\left[|\xi_i|\right]\right)^{p-1}\left(\sum_{i=1}^n \mathbb{E}\left[|\xi_i|\right]\right)\\
	&\le \frac{1}{p^{(p-1)\alpha}}\left(\sum_{i=1}^n \mathbb{E}\left[|\xi_i|\right]\right)^p.
	\end{align*}
	Therefore (in case 1),
	\begin{equation}\label{eq:Case1}
	p^{p\alpha}\sum_{i=1}^n \mathbb{E}\left[|\xi_i|^p\right] \le 2p^{p\alpha}\mathbb{E}\left[\max_{1\le i\le n}|\xi_i|^p\right] + 2p^{\alpha}\left(\sum_{i=1}^n \mathbb{E}\left[|\xi_i|\right]\right)^p.
	\end{equation}
	\item[--] \textit{Case 2:} If
	\[
	p^{p\alpha}\mathbb{E}\left[\max_{1\le i\le n}|\xi_i|^p\right] \ge \left(\sum_{i=1}^n \mathbb{E}\left[|\xi_i|\right]\right)^p,
	\]
	then
	\begin{align*}
	\mathbb{E}\left[\max_{1\le i\le n}|\xi_i|^{p-1}\right]&\left(\sum_{i=1}^n \mathbb{E}\left[|\xi_i|\right]\right)\\ &\le \left(\mathbb{E}\left[\max_{1\le i\le n}|\xi_i|^p\right]\right)^{(p-1)/p}\left(\sum_{i=1}^n \mathbb{E}\left[|\xi_i|\right]\right)\\
	&\le \left(\mathbb{E}\left[\max_{1\le i\le n}|\xi_i|^p\right]\right)^{(p-1)/p}p^{\alpha}\left(\mathbb{E}\left[\max_{1\le i\le n}|\xi_i|^p\right]\right)^{1/p}\\
	&\le p^{\alpha}\mathbb{E}\left[\max_{1\le i\le n}|\xi_i|^p\right].
	\end{align*}
	Therefore (in case 2),
	\begin{align}
	p^{p\alpha}\sum_{i=1}^n \mathbb{E}\left[|\xi_i|^p\right] &\le 2p^{p\alpha}\mathbb{E}\left[\max_{1\le i\le n}|\xi_i|^p\right] + 2p^{\alpha(p+1)}\mathbb{E}\left[\max_{1\le i\le n}|\xi_i|^p\right]\nonumber\\
	&\le 2p^{p\alpha}\mathbb{E}\left[\max_{1\le i\le n}|\xi_i|^p\right] + 2p^{p\alpha}\left(e^{1/e}\right)^{p\alpha}\mathbb{E}\left[\max_{1\le i\le n}|\xi_i|^p\right]\nonumber\\
	&\le (2 + (1.5)^{p\alpha})p^{p\alpha}\mathbb{E}\left[\max_{1\le i\le n}|\xi_i|^p\right]\label{eq:Case2}.
	\end{align}
\end{itemize}
Combining inequalities~\eqref{eq:Case1} and~\eqref{eq:Case2}, we get for $p\ge 1$ and $\alpha > 0$ that
\begin{align*}
p^{p\alpha}\sum_{i=1}^n \mathbb{E}\left[|\xi_i|^p\right] &\le (2 + (1.5)^{p\alpha})p^{p\alpha}\mathbb{E}\left[\max_{1\le i\le n}|\xi_i|^p\right] + 2p^{\alpha}\left(\sum_{i=1}^n \mathbb{E}\left[|\xi_i|\right]\right)^p\\
&\le 4(1.5)^{p\alpha}p^{p\alpha}\mathbb{E}\left[\max_{1\le i\le n}|\xi_i|^p\right] + 2p^{\alpha}\left(\sum_{i=1}^n \mathbb{E}\left[|\xi_i|\right]\right)^p\\
&\le 4(1.5)^{p\alpha}p^{p\alpha}\mathbb{E}\left[\max_{1\le i\le n}|\xi_i|^p\right] + 2(p^{1/p})^{p\alpha}\left(\sum_{i=1}^n \mathbb{E}\left[|\xi_i|\right]\right)^p\\
&\le 4(1.5)^{p\alpha}p^{p\alpha}\mathbb{E}\left[\max_{1\le i\le n}|\xi_i|^p\right] + 2(1.5)^{p\alpha}\left(\sum_{i=1}^n \mathbb{E}\left[|\xi_i|\right]\right)^p.
\end{align*}
This proves the result.
\end{proof}
\begin{lem}\label{lem:BoundingB}
Under the notation of Theorem~\ref{thm:MainUStat}, the quantity $B$ defined in~\eqref{eq:DefinitionsABCD} satisfies
\begin{align*}
B &\le \sup\left\{\sum_{1\le i\neq j\le n}\mathbb{E}\left[q_i(X_i)\sigma_{i,\phi}(X_i)w_{i,j}(X_i, X_j)\sigma_{j,\psi}(X_j)p_j(X_j)\right]:\right.\\
&\qquad\qquad\quad\left.\sum_{j = 1}^n \mathbb{E}\left[q_i^2(X_i)\right] \le 1,\,\sum_{i=1}^n \mathbb{E}\left[p_j^2(X_j)\right] \le 1\right\}.
\end{align*}
\end{lem}
\begin{proof}
Following the proof of Theorem 3.2 of~\cite{Gine00}, the quantity $B$ is the square root of the wimpy variance of
\[
S_n := \left(\sum_{i=1}^n \mathbb{E}\left[F_i^2(\varepsilon_i, Z_i; \mathcal{Z}_n')\big|\mathcal{Z}_n'\right]\right)^{1/2},
\]
where $\mathcal{Z}_n' := \{(\varepsilon_1', Z_1'), \ldots, (\varepsilon_n', Z_n')\}$ and
\[
F_i(\varepsilon_i, Z_i; \mathcal{Z}_n') := \varepsilon_i\Phi_{i,1}\sum_{j = 1, j\neq i}^n w_{i,j}(X_i, X_j')\Psi'_{j,1}\varepsilon_j'.
\]
This implies that
\[
S_n \le \left(\sum_{i=1}^n \mathbb{E}\left[G_i^2(X_i; \mathcal{Z}_n')\big|\mathcal{Z}_n'\right]\right)^{1/2},
\]
where for $\sigma_{i,\phi}^2(x) := \mathbb{E}\left[\phi^2(Y_i)\big|X_i = x\right]$,
\[
G_i(X_i; \mathcal{Z}_n') := \sigma_{i,\phi}(X_i)\sum_{j = 1, j\neq i}^n w_{i,j}(X_i, X_j')\Psi'_{j,1}\varepsilon_j'.
\]
Note that $\sigma_{i,\phi}(\cdot)$ depends on $i$ since the random variables are allowed to be non-identically distributed. Now observe that
\begin{equation}\label{eq:DualityFirst}
S_n = \sup\left\{\sum_{i=1}^n \int q_i(x)G_i(x; \mathcal{Z}_n')P_{X_i}(dx):\,\sum_{i=1}^n \int q_i^2(x)P_{X_i}(dx) \le 1\right\}.
\end{equation}
To prove this, note that for any $\{q_i(\cdot):\,1\le i\le n\}$ satisfying the (integral) constraint,
\begin{align*}
\sum_{i=1}^n &\int q_i(x)G_i(x; \mathcal{Z}_n')P_{X_i}(dx)\\
&\le \sum_{i=1}^n \left(\int q_i^2(x)P_{X_i}(dx)\right)^{1/2}\left(\int G_i^2(x; \mathcal{Z}_n')P_{X_i}(dx)\right)^{1/2}\\
&\le \left(\sum_{i=1}^n \int q_i^2(x)P_{X_i}(dx)\right)^{1/2}\left(\sum_{i=1}^n \int G_i^2(x; \mathcal{Z}_n')P_{X_i}(dx)\right)^{1/2} \le S_n.
\end{align*}
To prove the reverse inequality, define for $1\le i\le n$,
\[
q_i(x) := G_i(x; \mathcal{Z}_n')\left(\sum_{i=1}^n \int G_i^2(x; \mathcal{Z}_n')P_{X_i}(dx)\right)^{1/2}.
\]
It is clear that $\{q_i(\cdot):\,1\le i\le n\}$ satisfy the integral constraint in~\eqref{eq:DualityFirst} and
\[
\sum_{i=1}^n \int q_i(x)G_i(x; \mathcal{Z}_n')P_{X_i}(dx) = S_n.
\]
This completes the proof of~\eqref{eq:DualityFirst}. Rewriting the representation~\eqref{eq:DualityFirst}, we get
\[
S_n = \sup_{\sum_{i=1}^n \int q_i^2(x)P_{X_i}(dx) \le 1}\,\sum_{j = 1}^n \varepsilon_j'\Psi'_{j,1}\left(\sum_{i=1, i\neq j}^n \int q_i(x)\sigma_{i,\phi}(x)w_{i,j}(x, X_j')P_{X_i}(dx)\right).
\]
This representation shows that $S_n$ is indeed the supremum of an empirical process. The wimpy variance of this supremum is given by
\begin{align*}
\sup_{\{q_i(\cdot)\}}&\mbox{Var}\left(\sum_{j = 1}^n \varepsilon_j'\Psi'_{j,1}\left(\sum_{i=1, i\neq j}^n \int q_i(x)\sigma_{i,\phi}(x)w_{i,j}(x, X_j')P_{X_i}(dx)\right)\right)\\
&\le \sup_{\{q_i(\cdot)\}} \sum_{j = 1}^n \mathbb{E}\left[\sigma_{j,\psi}^2(X_j')\left(\sum_{i=1, i\neq j}^n \int q_i(x)\sigma_{i,\phi}(x)w_{i,j}(x, X_j')P_{X_i}(dx)\right)^2\right]\\
&= \sup_{\{q_i(\cdot)\}} \sum_{j = 1}^n \mathbb{E}\left[\sigma_{j,\psi}^2(X_j)\left(\sum_{i=1, i\neq j}^n \int q_i(x)\sigma_{i,\phi}(x)w_{i,j}(x, X_j)P_{X_i}(dx)\right)^2\right].
\end{align*}
Now a duality argument implies that
\begin{align*}
\sup_{\{q_i(\cdot)\}} &\left(\sum_{j = 1}^n \mathbb{E}\left[\sigma_{j,\psi}^2(X_j)\left(\sum_{i=1, i\neq j}^n \int q_i(x)\sigma_{i,\phi}(x)w_{i,j}(x, X_j)P_{X_i}(dx)\right)^2\right]\right)^{1/2}\\
&= \sup\left\{\sum_{1\le i\neq j\le n}\mathbb{E}\left[q_i(X_i)\sigma_{i,\phi}(X_i)w_{i,j}(X_i, X_j)\sigma_{j,\psi}(X_j)p_j(X_j)\right]:\right.\\
&\qquad\quad\left.\sum_{i=1}^n \mathbb{E}\left[p_j^2(X_j)\right] \le 1, \sum_{j = 1}^n \mathbb{E}\left[q_i^2(X_i)\right] \le 1\right\}.
\end{align*}
Thus the result follows.
\end{proof}

%%%%%%%%%%%%%%%%%%%%%%%%%%%%%%%%%%%%%%%%%%%%%
%%%%%%%%%%%%%%%%%%%%%%%%%%%%%%%%%%%%%%%%%%%%%

\section{Proofs of Results in Section~\ref{sec:UProc}}\label{appsec:UProc}

\subsection{Proof of Theorem~\ref{thm:MainUProc}}\label{appsec:MainUProc}
Similar to $\mathcal{U}_n^{(\ell)}, 1\le \ell\le 4$ defined in the proof of Theorem~\ref{thm:MainUStat}, we define
\begin{equation}\label{eq:SplitDegenerateUStats}
\begin{split}
\mathcal{U}_n^{(1)}(\mathcal{W}) &:= \sup_{w\in\mathcal{W}}\left|\sum_{1\le i\neq j\le n}\varepsilon_i\Phi_{i,1}w_{i,j}(X_i, X_j')\Psi'_{j,1}\varepsilon_j'\right|,\\
\mathcal{U}_n^{(2)}(\mathcal{W}) &:= \sup_{w\in\mathcal{W}}\left|\sum_{1\le i\neq j\le n}\varepsilon_i\Phi_{i,2}w_{i,j}(X_i, X_j')\Psi'_{j,1}\varepsilon_j'\right|,\\
\mathcal{U}_n^{(3)}(\mathcal{W}) &:= \sup_{w\in\mathcal{W}}\left|\sum_{1\le i\neq j\le n}\varepsilon_i\Phi_{i,1}w_{i,j}(X_i, X_j')\Psi'_{j,2}\varepsilon_j'\right|,\\
\mathcal{U}_n^{(4)}(\mathcal{W}) &:= \sup_{w\in\mathcal{W}}\left|\sum_{1\le i\neq j\le n}\varepsilon_i\Phi_{i,2}w_{i,j}(X_i, X_j')\Psi'_{j,2}\varepsilon_j'\right|.
\end{split}
\end{equation}
As in the proof of Theorem~\ref{thm:MainUStat}, we will control each of the terms separately in the following lemmas. All the lemmas below assume~\ref{assump:ConditionalTailAssumption} and~\ref{assump:BoundedKernelsUniform}.
\begin{lem}[Control of $\mathcal{U}_n^{(4)}(\mathcal{W})$]\label{lem:Un4Bound}
There exists a constant $K > 0$ (depending only on $\alpha, \beta$) such that for all $p\ge1$,
\[
\norm{\mathcal{U}_n^{(4)}(\mathcal{W})}_p \le K\Lambda_2(\mathcal{W})p^{1/\alpha^* + 1/\beta^*}.
\]
\end{lem}
\begin{proof}
Since $\norm{w_{i,j}}_{\infty} \le B_{\mathcal{W}}$ for all $w\in\mathcal{W}$, it follows that
\[
\mathcal{U}_n^{(4)}(\mathcal{W}) \le B_{\mathcal{W}}\sum_{1\le i\neq j\le n}|\Phi_{i,2}\Psi_{j,2}'| \le B_{\mathcal{W}}\left(\sum_{i=1}^n |\Phi_{i,2}|\right)\left(\sum_{j = 1}^n |\Psi_{j,2}'|\right).
\]
By definition
\begin{align*}
\mathbb{P}\left(\max_{1\le I\le n}\sum_{i=1}^I |\Phi_{i,2}| > 0\big|\mathcal{X}_n\right) &\le \mathbb{P}\left(\max_{1\le i\le n}|\phi(Z_i)| \ge T_{\phi}\big|\mathcal{X}_n\right) \le 1/8,\\
\mathbb{P}\left(\max_{1\le I\le n}\sum_{i=1}^I |\Psi_{i,2}'| > 0\big|\mathcal{X}_n'\right) &\le \mathbb{P}\left(\max_{1\le i\le n}|\psi(Z_i')| \ge T_{\psi}\big|\mathcal{X}_n'\right) \le 1/8.
\end{align*}
Hence by (6.8) of~\cite{LED91}, we get that
\begin{align*}
\mathbb{E}\left[\sum_{i=1}^n |\Phi_{i,2}|\big|\mathcal{X}_n\right] &\le C\mathbb{E}\left[\max_{1\le i\le n}|\phi(Z_i)|\big|\mathcal{X}_n\right]\\
\mathbb{E}\left[\sum_{i=1}^n |\Psi_{i,2}'|\big|\mathcal{X}_n'\right] &\le C\mathbb{E}\left[\max_{1\le i\le n}|\psi(Z_i')|\big|\mathcal{X}_n'\right],
\end{align*}
for some constant $C > 0$. Thus by applying Theorem 6.21 of~\cite{LED91} to $\sum \{\Phi_{i,1} - \mathbb{E}[\Phi_{i,1}|\mathcal{X}_n]\}$ and $\sum \{\Psi_{i,2}' - \mathbb{E}[\Psi_{i,2}'|\mathcal{X}_n']\}$, we get
\begin{align*}
\norm{\sum_{i=1}^n |\Phi_{i,2}|}_{\psi_{\alpha^*}|\mathcal{X}_n} &\le C\norm{\max_{1\le i\le n}|\phi(Z_i)|}_{\psi_{\alpha}|\mathcal{X}_n} \le CC_{\phi}(\log n)^{1/\alpha},\\
\norm{\sum_{i=1}^n |\Psi_{i,2}'|}_{\psi_{\beta^*}|\mathcal{X}_n} &\le C\norm{\max_{1\le i\le n}|\psi(Z_i')|}_{\psi_{\beta}|\mathcal{X}_n'} \le CC_{\psi}(\log n)^{1/\beta},\\
\end{align*}
Therefore, for all $p\ge1$,
\[
\norm{\mathcal{U}_n^{(4)}(\mathcal{W})}_p \le KB_{\mathcal{W}}C_{\phi}C_{\psi}(\log n)^{\alpha^{-1} + \beta^{-1}}p^{1/\alpha^* + 1/\beta^*} = K\Lambda_2(\mathcal{W})p^{1/\alpha^* + 1/\beta^*}.
\]
This completes the proof.
\end{proof}
The following lemma controls the moments of $\mathcal{U}_n^{(2)}(\mathcal{W})$ and $\mathcal{U}_n^{(3)}(\mathcal{W})$.
\begin{lem}[Control of $\mathcal{U}_n^{(2)}(\mathcal{W})$ and $\mathcal{U}_n^{(3)}(\mathcal{W})$]\label{lem:Un23Bound}
There exists a constant $K > 0$ (depending only on $\alpha, \beta$) such that for $p\ge 1$,
\begin{align*}
\norm{\mathcal{U}_n^{(2)}(\mathcal{W})}_p &\le Kp^{1/\alpha^*}\left[E_{n,2}(\mathcal{W}) + (\log n)^{1/2}\Sigma_{n,2}^{1/2}(\mathcal{W}) + (\log n)\Lambda_2(\mathcal{W})\right]\\
&\qquad+ Kp^{1/2 + 1/\alpha^*}\Sigma_{n,2}^{1/2}(\mathcal{W}) + Kp^{1 + 1/\alpha^*}\Lambda_2(\mathcal{W})\\
\norm{\mathcal{U}_n^{(3)}(\mathcal{W})}_p &\le Kp^{1/\beta^*}\left[E_{n,1}(\mathcal{W}) + (\log n)^{1/2}\Sigma_{n,1}^{1/2}(\mathcal{W}) + (\log n)\Lambda_2(\mathcal{W})\right]\\
&\qquad+ Kp^{1/2 + 1/\beta^*}\Sigma_{n,1}^{1/2}(\mathcal{W}) + Kp^{1 + 1/\beta^*}\Lambda_2(\mathcal{W}).
\end{align*}
\end{lem}
\begin{proof}
We will only prove the bound for $\mathcal{U}_n^{(2)}(\mathcal{W})$ and the proof for $\mathcal{U}_n^{(3)}(\mathcal{W})$ follows very similar arguments. Recall that
\[
\mathcal{U}_n^{(2)}(\mathcal{W}) = \sup_{w\in\mathcal{W}}\left|\sum_{i=1}^n \varepsilon_i\Phi_{i,2}g_i(X_i; \mathcal{Z}_n', w)\right|,\;\mbox{where}\; g_i(x; \mathcal{Z}_n', w) := \sum_{j = 1, j\neq i}^n \Psi_{j,1}'\varepsilon_j'w_{i,j}(x, X_j').
\]
Here again (6.8) of~\cite{LED91} applies and we get
\[
\norm{\mathcal{U}_n^{(2)}(\mathcal{W})}_{\psi_{\alpha^*}\big|\mathcal{X}_n, \mathcal{Z}_n'} \le KC_{\phi}(\log n)^{1/\alpha}\max_{1\le i\le n}\sup_{w\in\mathcal{W}}\left|\sum_{j = 1, j\neq i}^n \varepsilon_j'\Psi_{j,1}'w_{i,j}(X_i, X_j')\right|.
\]
By a similar calculation, we get
\[
\norm{\mathcal{U}_n^{(3)}(\mathcal{W})}_{\psi_{\beta^*}\big|\mathcal{X}_n', \mathcal{Z}_n} \le KC_{\psi}(\log n)^{1/\beta}\max_{1\le j\le n}\sup_{w\in\mathcal{W}}\left|\sum_{i = 1, i\neq j}^n \varepsilon_i\Phi_{i,1}w_{i,j}(X_i, X_j')\right|.
\]
Thus, for $p\ge 1$,
\begin{equation}\label{eq:FirstBoundU23Proc}
\begin{split}
\mathbb{E}\left[|\mathcal{U}_n^{(2)}(\mathcal{W})|^p\right] &\le K^pC_{\phi}^p(\log n)^{p/\alpha}p^{p/\alpha^*}\mathbb{E}\left[\max_{1\le i\le n}\sup_{w\in\mathcal{W}}\left|\sum_{j = 1, j\neq i}^n \varepsilon_j'\Psi_{j,1}'w_{i,j}(X_i, X_j')\right|^p\right],\\
\mathbb{E}\left[|\mathcal{U}_n^{(3)}(\mathcal{W})|^p\right] &\le K^pC_{\psi}^p(\log n)^{p/\beta}p^{p/\beta^*}\mathbb{E}\left[\max_{1\le j\le n}\sup_{w\in\mathcal{W}}\left|\sum_{i = 1, i\neq j}^n \varepsilon_i\Phi_{i,1}w_{i,j}(X_i, X_j')\right|^p\right].
\end{split}
\end{equation}
The right hand side quantities involve supremum of bounded empirical processes for which Talagrand's inequality applies; see proposition 3.1 of~\cite{Gine00}. Observe that for any $x\in\mathfrak{X}$,
\begin{align*}
\max_{1\le j\le n}\sup_{w\in\mathcal{W}}|\Psi_{j,1}'w_{i,j}(x, X_j')| &\le C_{\psi}(\log n)^{1/\beta}B_{\mathcal{W}},\\
\max_{1\le i\le n}\sup_{w\in\mathcal{W}}|\Phi_{i,1}w_{i,j}(X_i, x)| &\le C_{\phi}(\log n)^{1/\alpha}B_{\mathcal{W}}.
\end{align*}
By proposition 3.1 of~\cite{Gine00}, we obtain for any $x\in\mathfrak{X}$ and $p\ge1$,
\begin{align*}
&\mathbb{E}\left[\sup_{w\in\mathcal{W}}\,|g_i(x; \mathcal{Z}_n', w)|^p\right]%\\ &\qquad\qquad
\le K^p\left\{\bar{E}_{n,2}^p(\mathcal{W}) + p^{p/2}\bar{\Sigma}_{n,2}^{p/2}(\mathcal{W}) + p^pC_{\psi}^p(\log n)^{p/\beta}B_{\mathcal{W}}^p\right\},
\end{align*}
where $\bar{E}_{n,2}(\mathcal{W}) = C_{\phi}^{-1}E_{n,2}(\mathcal{W})/(\log n)^{1/\alpha}$ and $\bar{\Sigma}_{n,2}^{1/2}(\mathcal{W}) = C_{\phi}^{-1}\Sigma_{n,2}^{1/2}(\mathcal{W})/(\log n)^{1/\alpha}$. Therefore, by following the argument that lead to~\eqref{eq:Un2Bound}, we get that
\begin{align}
&\mathbb{E}\left[\max_{1\le i\le n}\sup_{w\in\mathcal{W}}|g_i(X_i; \mathcal{Z}_n', w)|^p\right]\label{eq:MaxgiBound}\\ &\qquad\le K^p\left[\bar{E}_{n,2}^p(\mathcal{W}) + p^{p/2}\bar{\Sigma}_{n,2}^{p/2}(\mathcal{W}) + p^pC_{\psi}^p(\log n)^{p/\beta}B_{\mathcal{W}}^p\right]\nonumber\\
&\qquad\quad+ K^p\left[(\log n)^{p/2}\bar{\Sigma}_{n,2}^{p/2}(\mathcal{W}) + (\log n)^pC_{\psi}^p(\log n)^{p/\beta}B_{\mathcal{W}}^p\right].\nonumber
\end{align}
Substituting this in~\eqref{eq:FirstBoundU23Proc}, we get
\begin{align*}
\mathbb{E}\left[|\mathcal{U}_n^{(2)}(\mathcal{W})|^p\right] &\le K^pp^{p/\alpha^*}\left[E_{n,2}^p(\mathcal{W}) + p^{p/2}\Sigma_{n,2}^{p/2}(\mathcal{W}) + p^p\Lambda_2^p(\mathcal{W})\right]\\
&\quad+ K^pp^{p/\alpha^*}\left[(\log n)^{p/2}\Sigma_{n,2}^{p/2}(\mathcal{W}) + (\log n)^{p}\Lambda_2^p(\mathcal{W})\right].
\end{align*}
By a similar calculation, we get
\begin{align*}
\mathbb{E}\left[|\mathcal{U}_n^{(3)}(\mathcal{W})|^p\right] &\le K^pp^{p/\beta^*}\left[E_{n,1}^p(\mathcal{W}) + p^{p/2}\Sigma_{n,1}^{p/2}(\mathcal{W}) + p^p\Lambda_2^p(\mathcal{W})\right]\\
&\quad+ K^pp^{p/\beta^*}\left[(\log n)^{p/2}\Sigma_{n,1}^{p/2}(\mathcal{W}) + (\log n)^{p}\Lambda_2^p(\mathcal{W})\right].
\end{align*}
This completes the proof of the result.
\end{proof}

The following lemma controls the moments of $\mathcal{U}_n^{(1)}(\mathcal{W})$. This is a bounded degenerate $U$-process and is (usually) the dominating term among the four parts.

\begin{lem}[Control of $\mathcal{U}_{n}^{(1)}(\mathcal{W})$]\label{lem:Un1Bound}
There exists a constant $K > 0$ (depending only on $\alpha, \beta$) such that for all $p\ge 1$,
\begin{align*}
\norm{\mathcal{U}_n^{(1)}(\mathcal{W})}_p
% &\le K^p\left(\mathbb{E}\left[\mathcal{U}_n^{(1)}(\mathcal{W})\right]\right)^p + K^pp^p\left[E_{n,2}^p(\mathcal{W}) + p^{p/2}\Sigma_{n,2}^{p/2}(\mathcal{W}) + p^p\Lambda_2^p(\mathcal{W})\right]\\
% &\quad+ K^pp^p\left[(\log n)^{p/2}\Sigma_{n,2}^{p/2}(\mathcal{W}) + (\log n)^p\Lambda_2^p(\mathcal{W})\right]\\
% &\quad+ K^pp^{p/2}\mathfrak{W}_{n,2}^p(\mathcal{W})\\
% &\quad+ K^pp^{p}\norm{(\phi w\psi)_{\mathcal{W}}}_{2\to2}^p + K^pp^{3p/2}\Sigma_{n,1}^{p/2}(\mathcal{W})\\
% &\quad+ K^pp^pE_{n,1}^p(\mathcal{W})\\
% &\quad+ K^pp^{p/2}\mathfrak{W}^p_{n,1}(\mathcal{W}) \\
&\le K\mathbb{E}\left[\mathcal{U}_n^{(1)}(\mathcal{W})\right] + Kp^{1/2}\left(\mathfrak{W}_{n,1}(\mathcal{W}) + \mathfrak{W}_{n,2}(\mathcal{W})\right)\\
&\quad+ Kp\left(\norm{(\phi w\psi)_{\mathcal{W}}}_{2\to2} + E_{n,1}(\mathcal{W}) + E_{n,2}(\mathcal{W}) + \Sigma_{n,2}^{1/2}(\mathcal{W})\sqrt{\log n} + \Lambda_2(\mathcal{W})\log n\right)\\
&\quad+ Kp^{3/2}\left(\Sigma_{n,1}^{1/2}(\mathcal{W}) + \Sigma_{n,2}^{1/2}(\mathcal{W})\right) + Kp^{2}\Lambda_{2}(\mathcal{W}).
\end{align*}
\end{lem}

\begin{proof}
Recall that
\[
\mathcal{U}_n^{(1)}(\mathcal{W}) = \sup_{w\in\mathcal{W}}\left|\sum_{i=1}^n \varepsilon_i\Phi_{i,1}g_i(X_i; \mathcal{Z}_n', w)\right|,\,\mbox{where}\, g_i(X_i; \mathcal{Z}_n', w) := \sum_{j = 1, j\neq i}^n \varepsilon_j'\Psi_{j,1}'w_{i,j}(X_i, X_j').
\]
Observe that conditional on $\mathcal{Z}_n'$, $\mathcal{U}_n^{(1)}(\mathcal{W})$ is a bounded empirical process and so Talagrand's inequality applies. Thus by Proposition 3.1 of~\cite{Gine00}, we get for $p\ge 1$
\begin{align*}
\mathbb{E}\left[|\mathcal{U}_n^{(1)}(\mathcal{W})|^p\big|\mathcal{Z}_n'\right] &\le K^p\left(\mathbb{E}\left[|\mathcal{U}_n^{(1)}(\mathcal{W})|\big|\mathcal{Z}_n'\right]\right)^p\\
&\quad+ K^pp^{p/2}\sup_{w\in\mathcal{W}}\left(\sum_{i=1}^n \mathbb{E}\left[\Phi_{i,1}^2g_i^2(X_i; \mathcal{Z}_n', w)\big|\mathcal{Z}_n'\right]\right)^{p/2}\\
&\quad+ K^pp^p\mathbb{E}\left[\max_{1\le i\le n}|\Phi_{i,1}|^p\sup_{w\in\mathcal{W}}\left|g_i(X_i; \mathcal{Z}_n', w)\right|^p\big|\mathcal{Z}_n'\right].
\end{align*}
Therefore, for $p\ge 1$,
\begin{align*}
\mathbb{E}\left[|\mathcal{U}_n^{(1)}(\mathcal{W})|^p\right] &\le K^p\mathbb{E}\left(\mathbb{E}\left[|\mathcal{U}_n^{(1)}(\mathcal{W})|\big|\mathcal{Z}_n'\right]\right)^p\\
&\quad+ K^pp^{p/2}\mathbb{E}\left[\sup_{w\in\mathcal{W}}\left(\sum_{i=1}^n \mathbb{E}\left[\sigma_{i,\phi}^2(X_i)g_i^2(X_i; \mathcal{Z}_n', w)\big|\mathcal{Z}_n'\right]\right)^{p/2}\right]\\
&\quad+ K^pp^pC_{\phi}^p(\log n)^{p/\alpha}\mathbb{E}\left[\max_{1\le i\le n}\sup_{w\in\mathcal{W}}\left|g_i(X_i; \mathcal{Z}_n', w)\right|^p\right]\\
&=: K^p\left[\mathbf{I} + \mathbf{II} + \mathbf{III}\right].
\end{align*}
\textit{Controlling $\mathbf{III}:$} Using~\eqref{eq:MaxgiBound} from Lemma~\ref{lem:Un23Bound}, we get
\begin{align*}
\mathbf{III} &\le K^pp^p\left[E_{n,2}^p(\mathcal{W}) + p^{p/2}\Sigma_{n,2}^{p/2}(\mathcal{W}) + p^p\Lambda_2^p(\mathcal{W})\right]\\
&\quad+ K^pp^p\left[(\log n)^{p/2}\Sigma_{n,2}^{p/2}(\mathcal{W}) + (\log n)^p\Lambda_2^p(\mathcal{W})\right].
\end{align*}
\textit{Controlling $\mathbf{II}:$} To control~$\mathbf{II}$, we use a technique similar to the one used in Lemma~\ref{lem:BoundingB}. For this note by~\eqref{eq:DualityFirst} that for any $w(\cdot, \cdot)$
\begin{align*}
&\left(\sum_{i=1}^n \int \sigma_{i,\phi}^2(x)g_i^2(x; \mathcal{Z}_n', w)P_{X_i}(dx)\right)^{1/2}\\ &\qquad= \sup\left\{\sum_{i=1}^n \int q_i(x)\sigma_{i,\phi}(x)g_i(x; \mathcal{Z}_n', w)P_{X_i}(dx):\,\sum_{i=1}^n \int q_i^2(x)P_{X_i}(dx) \le 1\right\}.
\end{align*}
Therefore,
\begin{align*}
\mathbf{II} &= p^{p/2}\mathbb{E}\left[\sup_{w\in\mathcal{W}}\sup_{\sum_{i=1}^n \int q_i^2(x)P_{X_i}(dx) \le 1}\left|\sum_{i=1}^n \int q_i(x)\sigma_{i,\phi}(x)g_i(x; \mathcal{Z}_n', w)P_{X_i}(dx)\right|^p\right].
\end{align*}
Now observe that
\begin{align*}
\sum_{i=1}^n \int q_i(x)\sigma_{i,\phi}(x)g_i(x; \mathcal{Z}_n', w)P_{X_i}(dx) = \sum_{j = 1}^n \varepsilon_j'\Psi_{j,1}'\ell_j(X_j'; \{q_i\}, w),
\end{align*}
where $\{q_i\}$ represents the sequence $(q_1, \ldots, q_n)$ satisfying $\sum_{i=1}^n \int q_i^2(x)P_{X_i}(dx) \le 1\}$ and
\begin{equation}\label{eq:elljDef}
\ell_j(X_j'; \{q_i\}, w) := \sum_{i=1, i\neq j}^n \int q_i(x)\sigma_{i,\phi}(x)w_{i,j}(x, X_j')P_{X_i}(dx).
\end{equation}
Thus
\[
\mathbf{II} = p^{p/2}\mathbb{E}\left[\sup_{w\in\mathcal{W}}\sup_{\{q_i\}}\left|\sum_{j = 1}^n \varepsilon_j'\Psi_{j,1}'\ell_j(X_j';\{q_i\}, w)\right|^p\right].
\]
The right hand side is a bounded empirical process and by proposition 3.1 of~\cite{Gine00}, we get
\begin{equation}\label{eq:BoundingII}
\begin{split}
&\mathbb{E}\left[\sup_{w\in\mathcal{W}}\sup_{\{q_i\}}\left|\sum_{j = 1}^n \varepsilon_j'\Psi_{j,1}'\ell_j(X_j';\{q_i\}, w)\right|^p\right]\\
&\qquad\le K^p\left(\mathbb{E}\left[\sup_{\{q_i\}}\sup_{w\in\mathcal{W}}\left|\sum_{j = 1}^n \varepsilon_j'\Psi_{j,1}'\ell_j(X_j';\{q_i\}, w)\right|\right]\right)^p\\
&\qquad\quad+ K^pp^{p/2}\sup_{\{q_i\}}\sup_{w\in\mathcal{W}}\left(\mbox{Var}\left(\sum_{j = 1}^n \varepsilon_j'\Psi_{j,1}'\ell_j(X_j'; \{q_i\}, w)\right)\right)^{p/2}\\
&\qquad\quad+ K^pp^p\mathbb{E}\left[\sup_{\{q_i\}}\sup_{w\in\mathcal{W}}\max_{1\le j\le n}|\Psi_{j,1}'|^p|\ell_j(X_j';\{q_i\}, w)|^p\right].
\end{split}
\end{equation}
We will now control each of the three terms appearing in~\eqref{eq:BoundingII}. Using the fact $|\Psi_{j,1}'| \le KC_{\phi}(\log n)^{1/\beta}$, we get
\[
\mathbb{E}\left[\sup_{\{q_i\}}\sup_{w\in\mathcal{W}}\max_{1\le j\le n}|\Psi_{j,1}'|^p|\ell_j(X_j';\{q_i\}, w)|^p\right] \le C_{\psi}^p(\log n)^{p/\beta}\sup_{w\in\mathcal{W}}\sup_{x'\in\mathfrak{X}}\sup_{\{q_i\}}|\ell_j(x'; \{q_i\}, w)|^p.
\]
By following the duality argument~\eqref{eq:DualityFirst}, we get
\[
\sup_{\{q_i\}}|\ell_j(x'; \{q_i\}, w)| \le \left(\sum_{i=1, i\neq j}^n \mathbb{E}\left[\sigma_{i,\phi}^2(X_i)w^2(X_i, x')\right]\right)^{1/2},
\]
and so,
\begin{equation}\label{eq:BoundingIIThirdTerm}
\begin{split}
&\mathbb{E}\left[\sup_{\{q_i\}}\sup_{w\in\mathcal{W}}\max_{1\le j\le n}|\Psi_{j,1}'|^p|\ell_j(X_j';\{q_i\}, w)|^p\right]\\ &\qquad\le K^pC_{\psi}^p(\log n)^{p/\beta}\sup_{w\in\mathcal{W}, x\in\mathfrak{X}}\,\left(\sum_{i=1}^n \mathbb{E}\left[\sigma_{i,\phi}^2(X_i)w^2(X_i, x')\right]\right)^{p/2} = K^p\Sigma_{n,1}^{p/2}(\mathcal{W}).
\end{split}
\end{equation}

Also, note that
\[
\mbox{Var}\left(\sum_{j = 1}^n \varepsilon_j'\Psi_{j,1}'\ell_j(X_j'; \{q_i\}, w)\right) = \sum_{j = 1}^n \mathbb{E}\left[\sigma_{j,\psi}^2(X_j')\ell_j^2(X_j';\{q_i\}, w)\right].
\]
Hence, again following the duality argument~\eqref{eq:DualityFirst}, we get
\begin{equation}\label{eq:BoundingIISecondTerm}
\begin{split}
&\sup_{\{q_i\}}\sup_{w\in\mathcal{W}}\left(\mbox{Var}\left(\sum_{j = 1}^n \varepsilon_j'\Psi_{j,1}'\ell_j(X_j'; \{q_i\}, w)\right)\right)^{p/2}\\
&\qquad\le \sup_{\{q_i\}}\sup_{w\in\mathcal{W}}\sup_{\{p_j\}}\left(\sum_{1\le i\neq j\le n}\mathbb{E}\left[q_i(X_i)\sigma_{i,\phi}(X_i)w_{i,j}(X_i, X_j')\sigma_{j,\psi}(X_j')p_j(X_j')\right]\right)^p.
\end{split}
\end{equation}
Here $\{p_j\}$ represents a sequence $(p_1, \ldots, p_n)$ satisfying $\sum_{j = 1}^n \int p_j^2(x)P_{X_j}(dx) \le 1$.
\medskip

Substituting~\eqref{eq:BoundingIISecondTerm} and~\eqref{eq:BoundingIIThirdTerm} in~\eqref{eq:BoundingII}, we get
\begin{equation}\label{eq:BoundII}
\begin{split}
\mathbf{II} &\le K^pp^{p/2}\left(\mathbb{E}\left[\sup_{\{q_i\}}\sup_{w\in\mathcal{W}}\left|\sum_{j = 1}^n \varepsilon_j'\Psi_{j,1}'\ell_j(X_j';\{q_i\}, w)\right|\right]\right)^p\\
&\qquad\quad+ K^pp^{p}\norm{(\phi w\psi)_{\mathcal{W}}}_{2\to2}^p + K^pp^{3p/2}\Sigma_{n,1}^{p/2}(\mathcal{W}).
\end{split}
\end{equation}
\textit{Controlling $\mathbf{I}:$} We use Lemma~\ref{lem:Lemma2Adam06} (a restatement of Lemma 2 of~\cite{Adamczak06}) to control $\mathbf{I}.$ In the notation of Lemma~\ref{lem:Lemma2Adam06}, take
\[
W_j = (\varepsilon_j', Z_j'),\,T = (Z_1, \ldots, Z_n, \varepsilon_1, \ldots, \varepsilon_n),
\]
and for $w\in\mathcal{W}$,
\[
f_j^w(W_j, T) = \sum_{i = 1, i\neq j}^n \varepsilon_i\Phi_{i,1}w_{i,j}(X_i, X_j')\Psi_{j,1}'\varepsilon_j'.
\]
This implies
\[
S = \mathbb{E}_T\left[\sup_{w\in\mathcal{W}}\left|\sum_{j = 1}^n \sum_{i = 1, i\neq j}^n \varepsilon_i\Phi_{i,1}w_{i,j}(X_i, X_j')\Psi_{j,1}'\varepsilon_j'\right|\right] = \mathbb{E}\left[\mathcal{U}_n^{(1)}(\mathcal{W})\big|\mathcal{Z}_n'\right].
\]
Observe that $\mathbb{E}[S] = \mathbb{E}[\mathcal{U}_n^{(1)}(\mathcal{W})]$. Thus we get for $p\ge 1$
\begin{align}
\mathbb{E}\left[S^p\right] &\le K^p\left(\mathbb{E}[S]\right)^p + K^pp^{p/2}\Upsilon^p\label{eq:ConditionExpectationBound}\\ &\qquad+ K^pp^p\mathbb{E}\left[\max_{1\le j\le n}\left(\mathbb{E}\left[\sup_{w\in\mathcal{W}}\left|\sum_{i = 1, i\neq j}^n \varepsilon_i\Phi_{i,1}w_{i,j}(X_i, X_j')\Psi_{j,1}'\varepsilon_j'\right|\big|\mathcal{Z}_n'\right]\right)^p\right],\nonumber
\end{align}
where
\begin{align*}
\Upsilon &:= \sup_{q\in\mathcal{Q}}\left(\sum_{j = 1}^n \mathbb{E}\left[\left(\sum_{w \in \mathcal{W}} \mathbb{E}_{T}\left[f_j^w(W_j, T)q_j(T)\right]\right)^2\right]\right)^{1/2},
\end{align*}
with $\mathcal{Q}$ defined in Lemma~\ref{lem:Lemma2Adam06}. We now simplify the last two terms on the right hand side of~\eqref{eq:ConditionExpectationBound}. First observe that for the third term
\begin{align*}
&\mathbb{E}\left[\sup_{w\in\mathcal{W}}\left|\sum_{i = 1, i\neq j}^n \varepsilon_i\Phi_{i,1}w_{i,j}(X_i, X_j')\Psi_{j,1}'\varepsilon_j'\right|\big|\mathcal{Z}_n'\right]\\ &\qquad\le KC_{\psi}(\log n)^{1/\beta}\sup_{x\in\mathfrak{X}}\mathbb{E}\left[\sup_{w\in\mathcal{W}}\left|\sum_{i = 1, i\neq j}^n \varepsilon_i\Phi_{i,1}w_{i,j}(X_i, x)\right|\right] = KE_{n,1}(\mathcal{W}).
\end{align*}
To control $\Upsilon$, observe that
\[
\sum_{w\in\mathcal{W}}\mathbb{E}_T\left[f_j^{w}(W_j, T)q_j(T)\right] = \varepsilon_j'\Psi_{j,1}'\sum_{w\in\mathcal{W}}\mathbb{E}_T\left[q_j(T)\sum_{i = 1, i\neq j}^n \varepsilon_i\Phi_{i,1}w_{i,j}(X_i, X_j')\right].
\]
So, using the definition of $\sigma_{j,\psi}^2(\cdot)$, we get
\begin{align*}
\Upsilon &= \sup_{q\in\mathcal{Q}}\left(\sum_{j = 1}^n \mathbb{E}\left[\sigma_{j,\psi}^2(X_j')\left(\sum_{w\in\mathcal{W}}\mathbb{E}_T\left[q_j(T)\sum_{i = 1, i\neq j}^n \varepsilon_i\Phi_{i,1}w_{i,j}(X_i, X_j')\right]\right)^2\right]\right)^{1/2}\\
&\overset{(a)}{=}\sup_{\{p_j\}, q\in\mathcal{Q}}\sum_{j = 1}^n \mathbb{E}\left[p_j(X_j')\sigma_{j,\psi}(X_j')\sum_{w\in\mathcal{W}}\mathbb{E}_T\left[q_j(T)\sum_{i = 1, i\neq j}^n \varepsilon_i\Phi_{i,1}w_{i,j}(X_i, X_j')\right]\right]\\
&\overset{(b)}{=} \sup_{\{p_j\}}\mathbb{E}\left[\sup_{w\in\mathcal{W}}\sum_{i = 1}^n \varepsilon_i\Phi_{i,1}\left(\sum_{j = 1, j\neq i}^n \int p_j(x)\sigma_{j,\psi}(x)w_{i,j}(X_i, x)P_{X_j}(dx)\right)\right].
\end{align*}
Equality~(a) above follows from the duality argument~\eqref{eq:DualityFirst} while equality~(b) follows from the argument given in Lemma~\ref{lem:Lemma2Adam06}.

\end{proof}

\subsection{Auxiliary Lemmas Used in Theorem~\ref{thm:MainUProc}}\label{app:AuxLemmaUProc}

The following lemma is a rewording of Lemma 2 of~\cite{Adamczak06}. For this result, define the class of functions
\[
\mathcal{Q} := \left\{q(\cdot) = (q_1(\cdot), q_2(\cdot), \ldots):\,\sum_{k = 1}^{\infty} |q_k(T)| = 1\quad\mbox{for all}\quad T\right\}.
\]
The domain of functions in $\mathcal{Q}$ is left out on purpose.
\begin{lem}\label{lem:Lemma2Adam06}
Suppose $\mathcal{F} := \{(f_1^k, \ldots, f_n^k):\,k\ge 1\}$ represents a countable class of vector functions. Define for independent random variables $T, W_1, \ldots, W_n$,
\[
S := \mathbb{E}_T\left[\sup_{k \ge 1}\left|\sum_{j = 1}^n f_j^k(W_j, T)\right|\right],
\]
where $\mathbb{E}_T[\cdot]$ represents the expectation only with respect to $T$. (So, $S$ is a random variable that depends on $W_1, \ldots, W_n$). If $\mathbb{E}_W[f_j^k(W_j, T)] = 0$ for a.e $T$, then there exists a constant $K > 0$ such that for all $p\ge 1$,
\begin{align*}
\mathbb{E}\left[S^p\right] &\le K^p(\mathbb{E}[S])^p + K^pp^{p/2}\sup_{q\in\mathcal{Q}}\left(\sum_{j = 1}^n \mathbb{E}\left[\left(\sum_{k = 1}^{\infty}\mathbb{E}_T[f_j^k(W_j, T)q_j(T)]\right)^2\right]\right)^{p/2}\\ &\qquad+ K^pp^p\mathbb{E}\left[\max_{1\le j\le n}\left(\mathbb{E}_T\left[\sup_{k\ge 1}|f_j^k(W_j, T)|\right]\right)^p\right].
\end{align*}
\end{lem}
\begin{proof}
Following the proof of Lemma 2 of~\cite{Adamczak06}, we get
\[
S = \sup_{q\in\mathcal{Q}}\,\left|\sum_{k = 1}^{\infty}\mathbb{E}_T\left[q_k(Y)\sum_{j = 1}^n f_j^k(W_j, T)\right]\right|.
\]
To see this, define $\hat{q}(\cdot) = (\hat{q}_1(\cdot), \ldots)\in\mathcal{Q}$ such that
\[
\hat{q}_{\hat{k}}(t) = \mbox{sign}\left(\sum_{j = 1}^n f_j^{\hat{k}}(W_j, T)\right),\quad\mbox{and}\quad \hat{q}_k(t) = 0,\quad\mbox{for $k\neq \hat{k}$.}
\]
Here $\hat{k}$ satisfying
\[
\left|\sum_{j = 1}^n f_j^{\hat{k}}(W_j, T)\right| = \sup_{k \ge 1}\left|\sum_{j = 1}^n f_j^{k}(W_j, T)\right|.
\]
Therefore,
\[
S = \sup_{q\in\mathcal{Q}}\left|\sum_{j = 1}^n \left(\sum_{k = 1}^{\infty} \mathbb{E}_T\left[q_k(T)f_j^k(W_j, T)\right]\right)\right| =: \sup_{q\in\mathcal{Q}}\left|\sum_{j = 1}^n g_{q,j}(W_j)\right|.
\]
The right hand side above is the supremum of a mean zero empirical process and so by proposition 3.1 of~\cite{Gine00}, we get
\begin{align*}
\mathbb{E}\left[S^p\right] &\le K^p(\mathbb{E}[S])^p + K^pp^{p/2}\sup_{q\in\mathcal{Q}}\left(\sum_{j = 1}^n \mathbb{E}\left[g_{q,j}^2(W_j)\right]\right)^{p/2} + K^pp^p\mathbb{E}\left[\max_{1\le j\le n}\sup_{q\in\mathcal{Q}}\left|g_{q,j}(W_j)\right|^p\right].
\end{align*}
From the definition of $\mathcal{Q}$, we get
\[
\sup_{q\in\mathcal{Q}}|g_{q,j}(W_j)| = \sup_{q\in\mathcal{Q}}\left|\sum_{k = 1}^{\infty}\mathbb{E}_T\left[q_k(T)f_j^k(W_j, T)\right]\right| = \mathbb{E}_T\left[\sup_{k\ge 1}|f_j^k(W_j, T)|\right].
\]
Thus,
\[
\mathbb{E}\left[\max_{1\le j\le n}\sup_{q\in\mathcal{Q}}\left|g_{q,j}(W_j)\right|^p\right] = \mathbb{E}\left[\max_{1\le j\le n}\left(\mathbb{E}_T\left[\sup_{k\ge 1}|f_j^k(W_j, T)|\right]\right)^p\right].
\]
So, the result follows.
\end{proof}

\section{Proof of the Maximal Inequality (Theorem~\ref{thm:UniformEntropyUProcesses})}\label{appsec:MaximalInequ}
The following moment bound of Rademacher chaos is used in the proof. See corollary 3.2.6 of \cite{DeLaPena99} and inequalities leading to (4.1.20) on page 167 of \cite{DeLaPena99}.
\begin{lem}\label{lem:Chaos}
Let $Z$ be a homogeneous Rademacher chaos of degree 2, that is,
\[
Z := \sum_{1\le i\neq j\le n} \epsilon_i\epsilon_ja_{i,j},
\]
for some constants $a_{i,j}, 1\le i\neq j\le n$. Then $\norm{Z}_{\psi_{1}} \le 4es_n$, where
\[
s_n^2 := \sum_{1\le i\neq j\le n} a_{i,j}^2.
\]
\end{lem}
\medskip
\begin{proof}[Proof of Theorem~\ref{thm:UniformEntropyUProcesses}]
As before, let $\mathcal{X}_n := \{X_1, X_2, \ldots, X_n\}.$ Also, let
\[
Z_{\epsilon}(f) := \left|\frac{1}{\sqrt{n(n-1)}}\sum_{1\le i\neq j\le n}\epsilon_i\epsilon_jf_{i,j}(X_i, X_j)\right|.
\]
By Lemma \ref{lem:Chaos}, we get conditional on $\mathcal{X}_n$,
\begin{align*}
\norm{Z_{\epsilon}(f)}_{\psi_{1}|\mathcal{Z}_n} &\le 4e\left(\frac{1}{n(n-1)}\sum_{1\le i\neq j\le n} f^2_{i,j}(X_i, X_j)\right)^{1/2} \le 4e\norm{f}_{2,P_n},
\end{align*}
where
\[
\norm{f}_{2, P_n} := \left(\frac{1}{n(n-1)}\sum_{1\le i\neq j\le n} f^2_{i,j}(X_i, X_j)\right)^{1/2},
\]
and define the discrete probability measure $P_n$ with support $\{X_1, \ldots, X_n\}$ as
\[
P_n(\{X_i\}) := \frac{1}{n}\quad\mbox{for}\quad 1\le i\le n.
\]
Now, following the proof of Theorem 5.1.4 of \cite{DeLaPena99},
\[
\norm{\max_{f\in\mathcal{F}} Z_{\epsilon}(f)}_{\psi_1|\mathcal{X}_n} \le C\int_0^{\Delta_n} \log N\left(\varepsilon, \mathcal{F}, \norm{\cdot}_{2,P_n}\right)d\varepsilon,
\]
where
\[
\Delta_n := \sup_{f\in\mathcal{F}}\norm{f}_{2,P_n}.
\]
Therefore,
\[
\norm{\max_{f\in\mathcal{F}} Z_{\epsilon}(f)}_{\psi_1|\mathcal{X}_n} \le C\norm{F}_{2,P_n}J_2\left(\frac{\Delta_n}{\norm{F}_{2,P_n}}, \mathcal{F}, \norm{\cdot}_{2}\right).
\]
This implies that
\begin{equation}\label{eq:ZeroBound}
\mathbb{E}\left[\sup_{f\in\mathcal{F}}Z_{\epsilon}(f)\right] \le C\mathbb{E}\left[\norm{F}_{2,P_n}J_2\left(\frac{\Delta_n}{\norm{F}_{2,P_n}}, \mathcal{F}, \norm{\cdot}_{2}\right)\right].
\end{equation}
Using concavity of $(x, y)\mapsto \sqrt{y}J_2(\sqrt{x/y}, \mathcal{F}, \norm{\cdot}_2)$ as in the proof of Theorem 2.1 of \cite{vdV11}, it follows that
\begin{equation}\label{eq:FirstBound}
\mathbb{E}\left[\norm{F}_{2,P_n}J_2\left(\frac{\Delta_n}{\norm{F}_{2,P_n}}, \mathcal{F}, \norm{\cdot}_2\right)\right] \le \norm{F}_{2,P}J_2\left(\frac{\sqrt{\mathbb{E}\left[\Delta_n^2\right]}}{\norm{F}_{2,P}}, \mathcal{F}, \norm{\cdot}_2\right),
\end{equation}
where
\[
\norm{F}_{2,P}^2 := \frac{1}{n(n-1)}\sum_{1\le i\neq j\le n}\mathbb{E}\left[F^2_{i,j}(X_i, X_j)\right].
\]
At this point the proof of Theorem 5.1 of \cite{CHEN17} uses Hoeffding averaging to bound $\mathbb{E}\left[\Delta_n^2\right]$ which proves the result for {\cred i.i.d.} %iid
random variables $X_i$. To allow for non-identically distributed random variables $X_i, 1\le i\le n$, we bound $\mathbb{E}[\Delta_n^2]$ in terms of $J_2$ on the right hand side of \eqref{eq:FirstBound}. This is similar to the proof of Theorem 2.1 of~\cite{vdV11}. To bound $\mathbb{E}\left[\Delta_n^2\right]$, define for $f\in\mathcal{F},$
\begin{align*}
W_n^{(1)}(f) &:= \frac{1}{n(n-1)}\left|\sum_{1\le i\neq j\le n} \left\{f^2_{i,j}(X_i, X_j) - \mathbb{E}\left[f^2_{i,j}(X_i, X_j)|X_i\right]\right\}\right.\\
&\qquad\qquad - \left.\vphantom{\sum_{1\le i\neq j\le n}}\left\{\mathbb{E}\left[f^2_{i,j}(X_i, X_j)|X_j\right] + \mathbb{E}\left[f^2_{i,j}(X_i, X_j)\right]\right\}\right|,\\
W_n^{(2)}(f) &:= \frac{1}{n(n-1)}\left|\sum_{1\le i\neq j\le n} \left\{\mathbb{E}\left[f^2_{i,j}(X_i, X_j)|X_i\right] - \mathbb{E}\left[f^2_{i,j}(X_i, X_j)\right]\right\}\right|,\\
W_n^{(3)}(f) &:= \frac{1}{n(n-1)}\left|\sum_{1\le i\neq j\le n}\left\{\mathbb{E}\left[f^2_{i,j}(X_i, X_j)|X_j\right] - \mathbb{E}\left[f^2_{i,j}(X_i, X_j)\right]\right\}\right|.
\end{align*}
Using these definitions, we get
\begin{equation}\label{eq:DeltanBound}
\Delta_n^2 \le \sup_{f\in\mathcal{F}}W_n^{(1)}(f) + \sup_{f\in\mathcal{F}}W_n^{(2)}(f) + \sup_{f\in\mathcal{F}} W_n^{(3)}(f) + \Sigma_n^2(\mathcal{F}),
\end{equation}
where
\[
\Sigma_n^2(\mathcal{F}) := \sup_{f\in\mathcal{F}}\frac{1}{n(n-1)}\sum_{1\le i\neq j\le n}\mathbb{E}\left[f^2_{i,j}(X_i, X_j)\right].
\]
By decoupling and symmetrization, we obtain
\begin{equation*}
\mathbb{E}\left[\sup_{f\in\mathcal{F}}W_n^{(1)}(f)\right] \le C\mathbb{E}\left[\sup_{f\in\mathcal{F}}\frac{1}{n(n-1)}\left|\sum_{1\le i\neq j\le n} \epsilon_i\epsilon_jf^2_{i,j}(X_i, X_j)\right|\right].
\end{equation*}
Set for $f\in\mathcal{F}$,
\[
R_{\epsilon}(f) := \frac{1}{\sqrt{n(n-1)}}\sum_{1\le i\neq j\le n} \epsilon_i\epsilon_jf^2_{i,j}(X_i, X_j).
\]
Again by Lemma \ref{lem:Chaos} and using $|f_{i,j}(x, x') + g_{i,j}(x, x')| \le 2R$ for all $f,g\in\mathcal{F}$ and $x, x'\in\mathcal{X}$, we get
\begin{align*}
\norm{R_{\epsilon}(f) - R_{\epsilon}(g)}_{\psi_1\big|\mathcal{X}_n} &\le 8eR\left(\frac{1}{n(n-1)}\sum_{1\le i\neq j\le n}\left(f_{i,j}(X_i, X_j) - g_{i,j}(X_i, X_j)\right)^2\right)^{1/2}\\
&\le 8eR\norm{f - g}_{2,P_n}.
\end{align*}
Hence by following the first part of the proof, we get
\begin{equation}\label{eq:W_n1}
\mathbb{E}\left[\sup_{f\in\mathcal{F}} W_n^{(1)}(f)\right] \le C\frac{R\norm{F}_{2,P}}{n}J_2\left(\frac{\sqrt{\mathbb{E}\left[\Delta_n^2\right]}}{\norm{F}_{2,P}}, \mathcal{F}, \norm{\cdot}_2\right).
\end{equation}
Substituting this in \eqref{eq:DeltanBound} after taking expectations,
\begin{align*}
\frac{\norm{\Delta_n}_2^2}{\norm{F}_{2,P}^2} \le CB^2_nJ_2\left(\frac{\norm{\Delta_n}_2}{\norm{F}_{2,P}}, \mathcal{F}, \norm{\cdot}_2\right) + A^2_n,
\end{align*}
where
\[
B^2_n := \frac{R}{n\norm{F}_{2,P}}\quad\mbox{and}\quad A^2_n := \frac{\mathbb{E}\left[\sup_{f\in\mathcal{F}}W_n^{(2)}(f)\right] + \mathbb{E}\left[\sup_{f\in\mathcal{F}}W_n^{(3)}(f)\right] + \Sigma_n^2(\mathcal{F})}{\norm{F}_{2,P}^2}.
\]
It follows that
\[
\frac{\norm{\Delta_n}_2^2}{\norm{F}_{2,P}^2} \le Cb^2_nJ_2\left(\frac{\norm{\Delta_n}_2}{\norm{F}_{2,P}}, \mathcal{F}, \norm{\cdot}_2\right) + a^2,
\]
for any $a \ge A_n$ and $b \ge B_n$. Therefore, by Lemma 2.1 of \cite{vdV11}, it follows that for any $a \ge A_n$ and $b \ge B_n$,
\[
J_2\left(\frac{\norm{\Delta_n}_2}{\norm{F}_{2,P}}, \mathcal{F}, \norm{\cdot}_2\right) \le CJ_2(a, \mathcal{F}, \norm{\cdot}_2)\left[1 + \frac{J_2(a, \mathcal{F}, \norm{\cdot}_2)b^2}{a^2}\right].
\]
Substituting this in \eqref{eq:FirstBound} and \eqref{eq:ZeroBound}, we get
\[
\mathbb{E}\left[\sup_{f\in\mathcal{F}} Z_{\epsilon}(f)\right] \le C\norm{F}_{2,P}J_2(a)\left[1 + \frac{J_2(a^2, \mathcal{F}, \norm{\cdot}_2)b^2}{a^2}\right],
\]
for any $a \ge A_n$ and $b\ge B_n$. The result is proved.% Note that the assumption all the functions $f\in\mathcal{F}$ are non-negative can be made without loss of generality since we can always split a real-valued function as a difference of non-negative functions.
\end{proof}

\bibliographystyle{apalike}
\bibliography{IntegralApprox}

\begin{thebibliography}{}

\bibitem[Adamczak, 2006]{Adamczak06}
Adamczak, R. (2006).
\newblock Moment inequalities for {$U$}-statistics.
\newblock {\em Ann. Probab.}, 34(6):2288--2314.

\bibitem[Adamczak and Kutek, 2023]{adamczak2023orlicz}
Adamczak, R. and Kutek, D. (2023).
\newblock On orlicz spaces satisfying the hoffmann-j $\{$$\backslash$o$\}$
  rgensen inequality.
\newblock {\em arXiv preprint arXiv:2310.04163}.

\bibitem[Arcones and Gin\'e, 1993]{Arcones93}
Arcones, M.~A. and Gin\'e, E. (1993).
\newblock Limit theorems for {$U$}-processes.
\newblock {\em Ann. Probab.}, 21(3):1494--1542.

\bibitem[Bakhshizadeh, 2023]{bakhshizadeh2023exponential}
Bakhshizadeh, M. (2023).
\newblock {Exponential tail bounds and large deviation principle for
  heavy-tailed $U$-statistics}.
\newblock {\em arXiv preprint arXiv:2301.11563}.

\bibitem[Bang and Robins, 2005]{bang2005doubly}
Bang, H. and Robins, J.~M. (2005).
\newblock Doubly robust estimation in missing data and causal inference models.
\newblock {\em Biometrics}, 61(4):962--973.

\bibitem[Bentkus and Götze, 1999]{BentkusGotze1999}
Bentkus, V. and Götze, F. (1999).
\newblock Optimal bounds in non-gaussian limit theorems for {$U$}-statistics.
\newblock {\em The Annals of Probability}, 27(1):454--521.

\bibitem[Boucheron et~al., 2005]{Bouch05}
Boucheron, S., Bousquet, O., Lugosi, G., and Massart, P. (2005).
\newblock Moment inequalities for functions of independent random variables.
\newblock {\em Ann. Probab.}, 33(2):514--560.

\bibitem[Boucheron et~al., 2013]{MR3185193}
Boucheron, S., Lugosi, G., and Massart, P. (2013).
\newblock {\em Concentration inequalities}.
\newblock Oxford University Press, Oxford.
\newblock A nonasymptotic theory of independence, With a foreword by Michel
  Ledoux.

\bibitem[Chakrabortty and Kuchibhotla, 2018]{chakrabortty2018tail}
Chakrabortty, A. and Kuchibhotla, A.~K. (2018).
\newblock Tail bounds for canonical {U}-statistics and {U}-processes with
  unbounded kernels.
\newblock Technical report, Working paper, Wharton School, University of
  Pennsylvania.

\bibitem[Chamakh et~al., 2021]{chamakh2021orlicz}
Chamakh, L., Gobet, E., and Liu, W. (2021).
\newblock Orlicz norms and concentration inequalities for $\beta$-heavy tailed
  random variables.
\newblock {\em Le Centre pour la Communication Scientifique Directe - HAL -
  Diderot}.

\bibitem[{Chen} and {Kato}, 2020]{CHEN17}
{Chen}, X. and {Kato}, K. (2020).
\newblock {Jackknife multiplier bootstrap: finite sample approximations to the
  $U$-process supremum with applications}.
\newblock {\em Probability Theory and Related Fields}, 176:1097--1163.

\bibitem[Cl{\'e}men\c{c}on et~al., 2008]{CLEM08}
Cl{\'e}men\c{c}on, S., Lugosi, G., and Vayatis, N. (2008).
\newblock Ranking and empirical minimization of {$U$}-statistics.
\newblock {\em Ann. Statist.}, 36(2):844--874.

\bibitem[de~la Pe{\~n}a, 1992]{DeLaPena92}
de~la Pe{\~n}a, V.~H. (1992).
\newblock Decoupling and {K}hintchine's inequalities for {$U$}-statistics.
\newblock {\em Ann. Probab.}, 20(4):1877--1892.

\bibitem[de~la Pe{\~n}a and Gin{\'e}, 1999]{DeLaPena99}
de~la Pe{\~n}a, V.~H. and Gin{\'e}, E. (1999).
\newblock {\em Decoupling}.
\newblock Probability and its Applications (New York). Springer-Verlag, New
  York.
\newblock From dependence to independence, Randomly stopped processes.
  $U$-statistics and processes. Martingales and beyond.

\bibitem[Delyon and Portier, 2016]{DEL16}
Delyon, B. and Portier, F. (2016).
\newblock Integral approximation by kernel smoothing.
\newblock {\em Bernoulli}, 22(4):2177--2208.

\bibitem[Dirksen, 2015]{Dirk15}
Dirksen, S. (2015).
\newblock Tail bounds via generic chaining.
\newblock {\em Electron. J. Probab.}, 20:no. 53, 29.

\bibitem[Gin{\'e} et~al., 2000]{Gine00}
Gin{\'e}, E., Lata{\l}a, R., and Zinn, J. (2000).
\newblock Exponential and moment inequalities for {$U$}-statistics.
\newblock In {\em High dimensional probability, {II} ({S}eattle, {WA}, 1999)},
  volume~47 of {\em Progr. Probab.}, pages 13--38. Birkh\"auser Boston, Boston,
  MA.

\bibitem[Gin{\'e} and Nickl, 2008]{Gine08}
Gin{\'e}, E. and Nickl, R. (2008).
\newblock A simple adaptive estimator of the integrated square of a density.
\newblock {\em Bernoulli}, 14(1):47--61.

\bibitem[Gin{\'e} and Nickl, 2016]{GINE16}
Gin{\'e}, E. and Nickl, R. (2016).
\newblock {\em Mathematical foundations of infinite-dimensional statistical
  models}.
\newblock Cambridge Series in Statistical and Probabilistic Mathematics.
  Cambridge University Press, New York.

\bibitem[Hall and Marron, 1987]{HallMarron87}
Hall, P. and Marron, J.~S. (1987).
\newblock Estimation of integrated squared density derivatives.
\newblock {\em Statist. Probab. Lett.}, 6(2):109--115.

\bibitem[He et~al., 2024]{he2024sparse}
He, Y., Wang, K., and Zhu, Y. (2024).
\newblock Sparse hanson-wright inequalities with applications.
\newblock {\em arXiv preprint arXiv:2410.15652}.

\bibitem[Houdr{\'e} and Reynaud-Bouret, 2003]{HOUDRE03}
Houdr{\'e}, C. and Reynaud-Bouret, P. (2003).
\newblock Exponential inequalities, with constants, for {U}-statistics of order
  two.
\newblock In {\em Stochastic inequalities and applications}, volume~56 of {\em
  Progr. Probab.}, pages 55--69. Birkh\"auser, Basel.

\bibitem[Kim, 2020]{kimstatistical2020}
Kim, I. (2020).
\newblock {\em Statistical Theory and Methods for Comparing Distributions}.
\newblock PhD thesis, Carnegie Mellon University.

\bibitem[Kolesko and Lata{\l}a, 2015]{Kole15}
Kolesko, K. and Lata{\l}a, R. (2015).
\newblock Moment estimates for chaoses generated by symmetric random variables
  with logarithmically convex tails.
\newblock {\em Statist. Probab. Lett.}, 107:210--214.

\bibitem[Kuchibhotla and Chakrabortty, 2022]{kuchibhotla2022moving}
Kuchibhotla, A.~K. and Chakrabortty, A. (2022).
\newblock {Moving beyond sub-Gaussianity in high-dimensional statistics:
  Applications in covariance estimation and linear regression}.
\newblock {\em Information and Inference: A Journal of the IMA},
  11(4):1389--1456.

\bibitem[Ledoux and Talagrand, 1991]{LED91}
Ledoux, M. and Talagrand, M. (1991).
\newblock {\em Probability in {B}anach spaces}, volume~23 of {\em Ergebnisse
  der Mathematik und ihrer Grenzgebiete (3) [Results in Mathematics and Related
  Areas (3)]}.
\newblock Springer-Verlag, Berlin.
\newblock Isoperimetry and processes.

\bibitem[Liu et~al., 2021]{liu2021adaptive}
Liu, L., Mukherjee, R., Robins, J.~M., and Tchetgen, E.~T. (2021).
\newblock Adaptive estimation of nonparametric functionals.
\newblock {\em Journal of Machine Learning Research}, 22(99):1--66.

\bibitem[Major, 2005]{Major05}
Major, P. (2005).
\newblock Tail behaviour of multiple random integrals and {$U$}-statistics.
\newblock {\em Probab. Surv.}, 2:448--505.

\bibitem[Major, 2013]{Major13}
Major, P. (2013).
\newblock {\em On the estimation of multiple random integrals and
  {$U$}-statistics}, volume 2079 of {\em Lecture Notes in Mathematics}.
\newblock Springer, Heidelberg.

\bibitem[Newey and Ruud, 2005]{NEW05}
Newey, W.~K. and Ruud, P.~A. (2005).
\newblock Density weighted linear least squares.
\newblock In {\em Identification and inference for econometric models}, pages
  554--573. Cambridge Univ. Press, Cambridge.

\bibitem[Nolan and Pollard, 1987]{Nolan87}
Nolan, D. and Pollard, D. (1987).
\newblock {$U$}-processes: rates of convergence.
\newblock {\em Ann. Statist.}, 15(2):780--799.

\bibitem[Nolan et~al., 1988]{nolan1988functional}
Nolan, D., Pollard, D., et~al. (1988).
\newblock Functional limit theorems for $ u $-processes.
\newblock {\em The Annals of Probability}, 16(3):1291--1298.

\bibitem[Robins et~al., 2008]{robins2008higher}
Robins, J., Li, L., Tchetgen, E., van~der Vaart, A., et~al. (2008).
\newblock Higher order influence functions and minimax estimation of nonlinear
  functionals.
\newblock In {\em Probability and statistics: essays in honor of David A.
  Freedman}, volume~2, pages 335--422. Institute of Mathematical Statistics.

\bibitem[Robins et~al., 2017]{robins2017minimax}
Robins, J.~M., Li, L., Mukherjee, R., Tchetgen, E.~T., and van~der Vaart, A.
  (2017).
\newblock Minimax estimation of a functional on a structured high-dimensional
  model.
\newblock {\em The Annals of Statistics}, 45(5):1951--1987.

\bibitem[Robins et~al., 2016]{robins2016asymptotic}
Robins, J.~M., Li, L., Tchetgen, E.~T., and van~der Vaart, A. (2016).
\newblock Asymptotic normality of quadratic estimators.
\newblock {\em Stochastic processes and their applications},
  126(12):3733--3759.

\bibitem[Robins et~al., 1994]{robins1994estimation}
Robins, J.~M., Rotnitzky, A., and Zhao, L.~P. (1994).
\newblock Estimation of regression coefficients when some regressors are not
  always observed.
\newblock {\em Journal of the American statistical Association},
  89(427):846--866.

\bibitem[Rudelson and Vershynin, 2013]{RudelsonVershynin13}
Rudelson, M. and Vershynin, R. (2013).
\newblock Hanson-{W}right inequality and sub-{G}aussian concentration.
\newblock {\em Electron. Commun. Probab.}, 18:no. 82, 9.

\bibitem[Serfling, 1980]{SERF80}
Serfling, R.~J. (1980).
\newblock {\em Approximation theorems of mathematical statistics}.
\newblock John Wiley \&\ Sons, Inc., New York.
\newblock Wiley Series in Probability and Mathematical Statistics.

\bibitem[Spokoiny and Zhilova, 2013]{Spokoiny13}
Spokoiny, V. and Zhilova, M. (2013).
\newblock Sharp deviation bounds for quadratic forms.
\newblock {\em Math. Methods Statist.}, 22(2):100--113.

\bibitem[Talagrand, 2014]{Tala14}
Talagrand, M. (2014).
\newblock {\em Upper and lower bounds for stochastic processes}, volume~60 of
  {\em Ergebnisse der Mathematik und ihrer Grenzgebiete. 3. Folge. A Series of
  Modern Surveys in Mathematics [Results in Mathematics and Related Areas. 3rd
  Series. A Series of Modern Surveys in Mathematics]}.
\newblock Springer, Heidelberg.
\newblock Modern methods and classical problems.

\bibitem[van~de Geer and Lederer, 2013]{Geer13}
van~de Geer, S. and Lederer, J. (2013).
\newblock The {B}ernstein-{O}rlicz norm and deviation inequalities.
\newblock {\em Probab. Theory Related Fields}, 157(1-2):225--250.

\bibitem[van~der Vaart and Wellner, 2011]{vdV11}
van~der Vaart, A. and Wellner, J.~A. (2011).
\newblock A local maximal inequality under uniform entropy.
\newblock {\em Electron. J. Stat.}, 5:192--203.

\bibitem[van~der Vaart and Wellner, 1996]{VdvW96}
van~der Vaart, A.~W. and Wellner, J.~A. (1996).
\newblock {\em Weak convergence and empirical processes}.
\newblock Springer Series in Statistics. Springer-Verlag, New York.
\newblock With applications to statistics.

\end{thebibliography}

\end{document}